%% file: spiralgibbs.tex
\documentclass[11pt]{article}
\input{macro.tex}

\title{Logarithm laws for equilibrium 
states\\ in negative curvature}
\author{Fr\'ed\'eric Paulin \and Mark Pollicott}

\date{\today} 

\begin{document} 
\bibliographystyle{../alphanum}
\maketitle

\begin{abstract} 
\noindent   
Let $M$ be a pinched negatively curved Riemannian manifold, whose unit
tangent bundle is endowed with a Gibbs measure $m_F$ associated to a
potential $F$. We compute the Hausdorff dimension of the conditional
measures of $m_F$. We study the $m_F$-almost sure asymptotic
penetration behaviour of locally geodesic lines of $M$ into small
neighbourhoods of closed geodesics, and of other compact (locally)
convex subsets of $M$.  We prove Khintchine-type and logarithm
law-type results for the spiraling of geodesic lines around these
objects. As an arithmetic consequence, we give almost sure Diophantine
approximation results of real numbers by quadratic irrationals with
respect to general Hölder quasi-invariant measures. \footnote{{\bf
    Keywords:} geodesic flow, negative curvature, spiraling,
  Khintchine theorem, logarithm law, equilibrium states, Gibbs
  measure, pressure, Hausdorff dimension, Diophantine approximation,
  quadratic irrational.  ~~{\bf AMS codes:} 37D35, 53D25, 37D40,
  53C22, 37A25, 37A45, 11J83, 28D20}
\end{abstract}

\section{Introduction} 
\label{sec:intro} 

Let $M$ be a complete connected Riemannian manifold with pinched
sectional curvature at most $-1$, and let $(\flow t)_{t\in\RR}$ be its
geodesic flow. In this paper, we consider for instance a closed
geodesic $D_0$ in $M$, and we want to study the spiraling of geodesics
lines around $D_0$. Given an ergodic probability measure $m$ invariant
under $(\flow t)_{t\in\RR}$, whose support is the nonwandering set
$\Omega$ of the geodesic flow, $m$-almost every orbit is dense in
$\Omega$. Two geodesic lines, having at some time their unit tangent
vectors very close, remain close for a long time. Hence $m$-almost
every geodesic line will stay for arbitrarily long periods of times in
a given small neighbourhood of $D_0$. In what follows, we make this
behaviour quantitative when $m$ is any equilibrium state.

Let $F:T^1M\ra\RR$ be a {\it potential}, that is, a
H\"older-continuous function.  Let $\M$ be the set of probability
measures $m$ on $T^1M$ invariant under the geodesic flow, for which
the negative part of $F$ is $m$-integrable, and let $h_m(\flow 1)$ be
the (metric) entropy of the geodesic flow with respect to $m$.
The {\it pressure} of the potential $F$ is
\begin{equation}\label{eq:defpressure}
P=P(F)=\sup_{m\in\M}\;\big(h_m(g^1)+\int_{T^1M}F\,dm\big)\;.
\end{equation} 
Let $m_F$ be a Gibbs measure on $T^1M$ associated to the potential $F$
(see \cite{PauPolSha} and Section \ref{sec:reminder}).  When finite
and normalised to be a probability measure (and if the negative part
of $F$ is $m_F$-integrable), it is the unique {\it equilibrium state},
that is, it attains the upper bound defining the pressure $P(F)$ (see
\cite[Theo.~6.1]{PauPolSha}, improving \cite{OtaPei04} when
$F=0$). For instance, $m_F$ is (up to a constant multiple) the
Bowen-Margulis measure $m_{\rm BM}$ if $F=0$, and is the Liouville
measure if $F$ is the strong unstable Jacobian $v\mapsto
-\,\frac{d}{dt}_{\mid t= 0}\log\operatorname{Jac} \big({\flow t}_{\mid
  W^{\rm su}(v)}\big)(v)$ and $M$ is compact (see
\cite[Theo.~7.2]{PauPolSha} for a more general case).  We will use the
construction of $m_F$ by Paulin-Pollicott-Schapira \cite{PauPolSha}
(building on work of Hamenstädt, Ledrappier, Coudène, Mohsen) via
Patterson densities $(\mu^F_x)_{x\in\wt M}$ on the boundary at
infinity $\partial_\infty\wt M$ of a universal cover $\wt M$ of $M$
associated to the potential $F$.

We first prove (see Section \ref{sec:ledrappier}) the following result
relating measure theoretic invariants of $m_F$ and $\mu^F_x$, which
extends Ledrappier's result \cite[\S 4]{Ledrappier13} when $F=0$.

\btheo \label{theo:ledrappintro} If $m_F$ is finite and $F$ is
$m_F$-integrable, the Hausdorff dimension of the Patterson measure
$\mu^F_x$ (with respect to the Gromov-Bourdon visual distance on
$\partial_\infty\wt M$) is equal to the metric entropy of the Gibbs
measure $m_F$ (for the geodesic flow).  
\etheo

Let $D_0$ be a closed geodesic in $M$ of length $\ell_0$. If $v_0\in
T^1M$ is tangent to $D_0$, let
\begin{equation}\label{eq:defPzero}
P_0=P(F_{\mid T^1D_0})=\frac{\max\{\int_0^{\ell_0}F(\flow t v_0)\,dt,
\int_0^{\ell_0}F(\flow t(- v_0))\,dt\}}{\ell_0}\;.
\end{equation} 
We will prove that $P_0<P$ if $m_F$ is finite. Let $\epsilon_0>0$ and
let $\psi:[0,+\infty\mathclose{[} \ra\mathopen{[}0,+\infty[$ be a
Lipschitz map. As introduced in \cite{HerPau10}, let $E(\psi)$ be the
set of {\it $(\epsilon_0,\psi)$-Liouville vectors around $D_0$}, that
is the set of $v\in T^1M$ such that there exists a sequence
$(t_n)_{n\in\NN}$ in $[0,+\infty[$ converging to $+\infty$ such that
for every $t\in [t_n, t_n + \psi(t_n)]$, the footpoint of $\flow t v$
belongs to the $\epsilon_0$-neighbourhood $\N_{\epsilon_0}D_0$ of $0$.

The Khintchine-type result describing the spiraling around the closed
geodesic $D_0$ is the following (simplified version of the) main
result of this paper (see Section \ref{sec:spiral}).

\btheo\label{theo:mainintro} Assume that $M$ is compact. If the
integral $\int_0^{+\infty} e^{\psi(t)(P_0-P)} \;dt$ diverges
(resp.~converges) then $m_F$-almost every (resp.~no) point of $T^1M$
belongs to $E(\psi)$.  
\etheo

When $F=0$ (that is, when $m_F$ is the Bowen-Margulis measure), this
result is due to Hersonsky-Paulin \cite{HerPau10}. As $m_F$ can be
taken to be the Liouville measure, this result answers a question
raised in loc.~cit. This result, in this particular case when $D_0$ is
a closed geodesic, can be restated as a well-approximation type of
result of points in the limit set of the fundamental group of $\Ga$ by
an orbit of a loxodromic fixed point, see for instance
\cite{FisSimUrb14} for very general results (their measure on the
limit set corresponds to $F=0$, though an extension might be
possible), and the references of \cite{FisSimUrb14} for historical
motivation and partial results. This result is a shrinking target
problem type, and our main tool is the mixing property of the geodesic
flow of $M$ for Gibbs measure (see \cite{PauPolSha}).

We stated this result as such to emphasize its novelty even in the
compact case, but it is true in a much more general setting, both from
$M$ and $D_0$ (see Theorem \ref{theo:maintechnicalspiralGibbs}).
For instance, when $M$ is a geometrically finite locally symmetric
orbifold, when $F$ has finite pressure $P(F)$ and finite Gibbs measure
$m_F$, when $D_0$ is a compact totally geodesic suborbifold (of
positive dimension and codimension), the result still holds. When $M$
is the quotient of real hyperbolic $3$-space by a geometrically
finite Kleinian group $\Ga$, when $F$ has finite pressure $P(F)$ and
finite Gibbs measure $m_F$, and when $D_0$ is the convex hull of the
limit set of a precisely invariant quasi-fuschian closed surface
subgroup $\Ga_0$ of $\Ga$, the result still holds.  See Section
\ref{sec:spiral} for more examples.

\medskip When $F=0$, the following logarithm law for the almost sure
spiraling of geodesic lines around $D_0$ is due to Hersonsky-Paulin
\cite{HerPau10}.  Let $\pi:T^1 M\ra M$ be the unit tangent
bundle. Define the penetration map $\ppp:T^1M\times \RR\ra[0,+\infty]$
of the geodesic lines inside $\N_{\epsilon_0}D_0$ by $\ppp(v,t)=0$ if
$\pi(\phi_tv)\notin \N_{\epsilon_0}D_0$, and otherwise $\ppp(v,t)$ is
the maximal length of an interval $I$ in $\RR$ containing $t$ such
that $\pi(\phi_sv)\in \N_{\epsilon_0}D_0$ for every $s\in I$.

\bcoro \label{coro:mainintro} Under the assumptions of Theorem
\ref{theo:mainintro}, for $m_F$-almost every $v\in T^1M$, we have
$$ 
\limsup_{t\ra+\infty}\frac{\ppp(v,t)}{\ln t} =
\frac{1}{P-P_0}\;.
$$
\ecoro

In Section \ref{sec:arithapplications}, we will give arithmetic
applications of Theorem \ref{theo:mainintro}. We will in particular
generalise to a huge class of measures on $\RR$ the Khintchine-type
result of approximation of real numbers by quadratic irrationals over
$\QQ$, that was proved in \cite{ParPau11MZ} for the Lebesgue measure,
and prove other $0$-$1$-laws of approximations of real numbers by
arithmetically defined points. To conclude this introduction, we give
one example of such a result.

Let $a,b\in\NN-\{0\}$ be positive integers such that the
equation $x^2-a\,y^2-b\,z^2=0$ has no nonzero integer solution (for
instance $a=2$ and $b=3$). Let $\Ga_{a,\,b}$ be
$$
\Big\{\begin{pmatrix} x+y\sqrt{a} & z-t\sqrt{a}\\
 b(z+t\sqrt{a}) & x-y\sqrt{a} 
\end{pmatrix}\,:\,(x,y,z,t)\in\ZZ^4\;\;{\rm and}\;\;
x^2-a\,y^2-b\,z^2+ab\,t^2=1\Big\}\;,
$$
which is a discrete subgroup of $\operatorname{SL}_2(\RR)$, whose
action by homographies on $\PP_1(\RR)=\RR\cup\{\infty\}$ is denoted by
$\cdot$~. If $\alpha\in\RR$ is a solution of the equation $\ga\cdot
X=X$ for some $\ga\in \Ga_{a,\,b}$, then $\alpha$ is quadratic over
$\QQ(\sqrt{a})$, and if furthermore $\alpha\notin \QQ(\sqrt{a})$, we
denote by $\alpha^\sigma$ its Galois conjugate over
$\QQ(\sqrt{a})$. Given $\ga\in \Ga_{a,\,b}$ with trace
$\operatorname{tr}\ga\neq 0,\pm 2$, we denote by $\ga^+$ and $\ga^-$
the attractive and repulsive fixed points of $\ga$ in
$\RR\cup\{\infty\}$.

Given a continuous action of a discrete group $G$ on a compact metric
space $(X,d)$, recall that a {\it Hölder quasi-invariant measure} (see
for instance \cite{Ledrappier95, Hamenstadt97}) on $X$ (for the action
of $G$) is a probability measure $\mu$ such that for every $g\in G$,
the measure $g_*\mu$ is absolutely continuous with respect to $\mu$,
and the Radon-Nykodim derivative $\frac{d\,g_*\mu}{d\,\mu}$ coincides
$\mu$-almost everywhere with a Hölder-continuous map on $X$, that we
will still denote by $\frac{d\,g_*\mu}{d\,\mu}$.

The next result is a Khintchine-type of result, under a huge class of
measures, for the Diophantine approximation of real numbers by
quadratic irrationals over $\QQ(\sqrt{a})$ in a (dense) orbit under
the arithmetic group $\Ga_{a,\,b}$ (extended by the Galois
conjugation).

\bcoro \label{coro:approxintro} Let $\mu$ be a Hölder quasi-invariant
measure on $\RR\cup\{\infty\}$ for the action by homographies of
$\Ga_{a,\,b}$.  Let $\ga_0$ be a primitive element in $\Ga_{a,\,b}$
with $\operatorname{tr}(\ga_0)\neq 0,\pm2$.  For $\mu$-almost every
$x\in\RR$, we have
$$
\;\liminf_{\alpha\in\Ga_{a,\,b}\cdot \{\ga_0^-,\ga_0^+\}\;:\;|\alpha-\alpha^\sigma|\ra 0}\;\;
\frac{|x-\alpha|}{|\alpha-\alpha^\sigma|(-\ln|\alpha-\alpha^\sigma|)^{-s}}=
0\;\;({\rm resp.} = +\infty)
$$ 
if $s\leq \frac1{\delta-\delta_0}$ (resp.~$s>
\frac1{\delta-\delta_0}$), where
$$
\delta= \limsup_{s\ra+\infty}\;\frac1{2\ln s} \;
\ln\sum_{\ga\in\Ga_{a,\,b},\;2<|\operatorname{tr}(\ga)|\leq s}
\frac{d(\ga^{-1})_*\mu}{d\,\mu}(\ga^+)
$$ 
and $\delta_0=
\frac1{2\operatorname{arcosh}(\frac{|\operatorname{tr}\ga_0|}{2})}\,
\max\big\{\frac{d(\ga_0^{-1})_*\mu}{d\,\mu}(\ga_0^+),\,
\frac{d(\ga_0)_*\mu}{d\,\mu}(\ga_0^-) \big\}$.  
\ecoro

We refer to Section \ref{sec:arithapplications} for more general
results, in particular for approximations with congruence properties
and for the approximation of complex numbers by quadratic irrationals
over an imaginary quadratic extension of $\QQ$.

\bigskip {\small {\it Acknowledgments: } The second author thanks the
  \'Ecole Normale Supérieure, for an Invited Professor position in
  2009 where this work was started, and the Université Paris-Sud
  (Orsay) for an Invited Professor position in 2014 where this work was
  completed.}

\section{A summary of the Patterson-Sullivan theory 
for Gibbs states}
\label{sec:reminder}

Most of the content of this section is extracted from
\cite{PauPolSha}, to which we refer for the proofs of the claims and
for more details.

Let $\wt M$ be a complete simply connected Riemannian manifold with
(dimension at least $2$ and) pinched negative sectional curvature
$-b^2\le K\le -1$, and let $x_0\in\wt M$ be a fixed basepoint. For
every $\epsilon>0$ and every subset $A$ of $\wt M$, we denote by
$\N_\epsilon A$ the closed $\epsilon$-neighbourhood of $A$ in $\wt M$.

We denote by $\pi:T^1 \wt M\ra \wt M$ the unit tangent bundle of $\wt
M$, where $T^1 \wt M$ is endowed with Sasaki's Riemannian metric. Let
$\partial_{\infty}\wt M$ be the boundary at infinity of $\wt M$. We
denote by $\Lambda G$ the limit set of any discrete group of
isometries $G$ of $\wt M$, and by $\C\Lambda G$ the convex hull in
$\wt M$ of $\Lambda G$, if $\Lambda G$ has at least two elements.

Let $\Ga$ be a nonelementary (not virtually nilpotent) discrete group
of isometries of $\wt M$. Let $M$ and $T^1M$ be the quotient
Riemannian orbifolds $\Ga\bs\wt M$ and $\Ga\bs T^1\wt M$, and let
again $\pi:T^1 M\ra M$ be the map induced by $\pi:T^1 \wt M\ra \wt
M$. We denote by $(\flow{t})_{t\in\RR}$ the geodesic flow on $T^1\wt
M$, as well as its quotient flow on $T^1M$.

For every $v\in T^1\wt M$, let $v_-\in\partial_\infty\wt M$ and
$v_+\in\partial_\infty\wt M$, respectively, be the endpoints at
$-\infty$ and $+\infty$ of the geodesic line $g_v:\RR \ra \wt M$
defined by $v$ (that is, such that $\dot{g_v}(0)=v$).  Let
$\partial_\infty^2\wt M$ be the subset of $\partial_\infty\wt
M\times\partial_\infty\wt M$ which consists of pairs of distinct
points at infinity of $\wt M$.  {\em Hopf's parametrisation} of
$T^1\wt M$ is the homeomorphism which identifies $T^1\wt M$ with
$\partial_\infty^2\wt M\times\RR$, by the map $v\mapsto(v_-,v_+,t)$,
where $t$ is the signed distance of the closest point to $x_0$ on
$g_v(\RR)$ to $\pi(v)$. Let $\wt\Omega\Ga$ be the $\Ga$-invariant set
of $v\in T^1\wt M$ such that $v_-,v_+\in\Lambda\Ga$, whose image in
$T^1M$ is the nonwandering set of the geodesic flow of $T^1M$.

Let $\iota: T^1\wt M\to T^1\wt M$ be the (Hölder-continuous) {\em
  antipodal (flip) map} of $T^1\wt M$ defined by $\iota v=-v$. We
denote the quotient map of $\iota$ again by $\iota: T^1 M\to T^1 M$.

Let $\wt F: T^1\wt M\to\RR$ be a fixed H\"older-continuous
$\Ga$-invariant function, called a {\em potential} on $T^1\wt M$. It
induces a H\"older-continuous function $F: T^1M \ra\RR$, called a {\it
  potential} on $T^1 M$.  Two potentials $\wt F$ and $\wt F^*$ on
$T^1\wt M$ (or their induced maps on $T^1M$) are {\it cohomologous} if
there exists a H\"older-continuous $\Ga$-invariant map $\wt G :T^1\wt
M\ra \RR$, differentiable along every flow line, such that
\begin{equation}\label{eq:cohomologouspotential}
\wt F^*(v)-\wt F(v)=\frac{d}{dt}_{\mid t=0}\wt G(\phi_tv)\;.
\end{equation}

For any two distinct points $x,y\in\wt M$, let $v_{xy}\in T^1_x \wt M$
be the initial tangent vector of the oriented geodesic segment $[x,y]$
in $\wt M$ that connects $x$ to $y$; define
$$
\int_x^y \wt F= \int_0^{d(x,y)}\wt F(\flow{t}v_{xy})\;dt\,
$$ 
and $\int_x^x \wt F=0$ for all $x\in\wt M$.  Given a hyperbolic
element $\ga\in\Ga$ with translation axis $A_\ga$, the {\it period} of
$\ga$ for $F$ is, for any $x\in A_\ga$,
$$
\operatorname{Per}_F(\ga)=\int_{x}^{\ga x} \wt F\;.
$$

The {\em critical
  exponent} of $(\Ga,F)$ is
$$
\delta_{\Ga,\,F}=\limsup_{n\to+\infty}\frac
1n\log\sum_{\ga\in\Ga, \, n-1<d(x,\ga y)\le n} e^{\int_x^{\ga y}\wt F}\;.
$$
When $F=0$, then $\delta_{\Ga,\,F}$ is the standard critical exponent
$\delta_\Ga$ of $\Ga$. Note that
\begin{equation}\label{symintFiota}
\int_x^y \wt F=\int_y^x \wt F\circ\iota\;,
\end{equation}
that $\delta_{\Ga,\,F}= \delta_{\Ga,\,F\circ\iota}>-\infty$ and that
$\delta_{\Ga,\,F+c}= \delta_{\Ga,\,F}+c$ for every constant $c>0$ (see
\cite[Lem.~3.3]{PauPolSha}). We assume that $\delta_{\Ga,\,F}<+\infty$
(this is for instance satisfied if $\wt F$ is bounded, see
\cite[Lem.~3.3]{PauPolSha}).  By \cite[Theo.~6.1]{PauPolSha}, the
critical exponent $\delta_{\Ga,\,F}$ is equal to the pressure $P(F)$
of $F$ on $T^1M$, defined in Equation \eqref{eq:defpressure}. The {\it
  Poincaré series}
$$ 
Q_{\Ga,\,F,\,x_0}(s)=\sum_{\ga\in\Ga}\;e^{\int_x^{\ga y}(\wt F-s)}
$$
of $(\Ga,F)$ converges if $s>\delta_{\Ga,\,F}$ and diverges if
$s<\delta_{\Ga,\,F}$. We say that $(\Ga, F)$ is {\it of divergence
  type} if $Q_{\Ga,\,F,\,x_0}(s)$ diverges at $s=\delta_{\Ga,\,F}$.

The {\em (normalised) Gibbs cocycle} of $\wt F$ is the function
$C^F:\partial_\infty\wt M\times\wt M\times\wt M\to\RR$ defined by
$$
(\xi,x,y)\mapsto C^F_\xi(x,y)
  =\lim_{t\to+\infty}\int_y^{\xi_t}(\wt F-\delta_{\Ga,\,F})-
\int_x^{\xi_t}(\wt F-\delta_{\Ga,\,F})\;, 
$$
where $t\mapsto\xi_t$ is any geodesic ray with endpoint
$\xi\in\partial_\infty\wt M$. We have $C^{F+c}=C^F$ for every constant
$c\in\RR$. If $\wt F=0$, then $C^F=\delta_{\Ga} \beta$, where $\beta$ is
the Busemann cocycle.

By \cite[Lem.~3.2, 3.4]{PauPolSha}, there exists a constant $c_1>0$
(depending only on the H\"older constants of $\wt F$ and on the bounds
of the sectional curvature of $\wt M$) such that for all $x,y,z\in \wt
M$, we have
\begin{equation}\label{eq:holdercontpotential}
\Big|\;\int_{x}^{z}\wt F- \int_{y}^{z}\wt F\;\Big|\leq 
c_1\;e^{d(x,\,y)}+ d(x,y)\max_{\pi^{-1}(B(x,\,d(x,\,y)))}|\wt F|\;,
\end{equation}
and, for every $\xi\in\partial_\infty\wt M$,
\begin{equation}\label{eq:holdercontcocycle}
\big|\;C^F_\xi(x,y)\;\big|\leq 
c_1\;e^{d(x,\,y)}+ d(x,y)\max_{\pi^{-1}(B(x,\,d(x,\,y)))}|\wt F|\;.
\end{equation}
\medskip

A family $(\mu^F_x)_{x\in\wt M}$ of finite measures on
$\partial_\infty \wt M$, whose support is the limit set $\Lambda\Ga$
of $\Ga$, is a {\em  Patterson density} for the pair
$(\Ga,\wt F)$ if
$$
\ga_*\mu^F_x=\mu^F_{\ga x}
$$
for all $\ga\in\Ga$ and $x\in \wt M$, and if the following
Radon-Nikodym derivatives exist for all $x,y\in\wt M$ and satisfy for
almost all $\xi\in\partial_\infty\wt M$
$$
\frac{d\mu^F_x}{d\mu^F_y}(\xi)=e^{-C^F_\xi(x,\,y)}\,.
$$

A {\em Gibbs measure} on $T^1\wt M$ for $(\Ga,\wt F)$ is the measure
$\wt m_{F}$ on $T^1\wt M$ given by the density
\begin{equation}\label{eq:defigibbs}
d\wt m_{F}(v)=
e^{C^{F\circ\iota}_{v_-}(x_0,\,\pi(v))\,+\,C^F_{v_+}(x_0,\,\pi(v))}\;
d\mu^{F\circ\iota}_{x_0}(v_-)\,d\mu^F_{x_0}(v_+)\,dt
\end{equation} 
in Hopf's parametrisation. Patterson densities $(\mu^F_x)_{x\in\wt M}$
and $(\mu^{F\circ\iota}_x)_{x\in\wt M}$ exist (see \cite[\S
3.6]{PauPolSha}, their construction, whence the existence of $\wt
m_F$, requires only $\Ga$ to be nonelementary and $\delta_{\Ga,\,F}
<+\infty$).  The Gibbs measure $\wt m_{F}$ is independent of $x_0$,
its support is $\wt \Omega\Ga$, and it is invariant under the actions
of the group $\Ga$ and of the geodesic flow.  Thus (see \cite[\S
2.6]{PauPolSha}), it defines a measure $m_{F}$ on $T^1M$ which is
invariant under the quotient geodesic flow, called a {\em Gibbs
  measure} on $T^1M$.  For every constant $c>0$, note that
$(\mu^F_x)_{x\in\wt M}$ is also a Patterson density for the pair
$(\Ga,\wt F+c)$, thus $\wt m_{F}$ is also a Gibbs measure for
$(\Ga,\wt F+c)$. If $m_{F}$ is finite, then the Patterson densities
are unique up to a common multiplicative constant (see \cite[\S
5.3]{PauPolSha}); hence the Gibbs measure of $m_F$ is uniquely
defined, up to a multiplicative constant, and, when normalised to be a
probability measure, it is the unique equilibrium state for the
potential $F$, if the negative part of $F$ is $m_F$-integrable, see
\cite[Theo.~6.1]{PauPolSha}.

By its definition as a quasi-product, the Gibbs measure $\wt m_F$
satisfies the following property, used without mention in what follows:
for every $x\in \wt M$, the premiage by $v\mapsto v_+$ of a set of
measure $0$ (respectively $>0$) for $\mu^F_x$ has measure $0$
(respectively $>0$) for $\wt m_F$.

We refer to \cite[\S 8]{PauPolSha} for finiteness criteria of $m_F$,
in particular satisfied if $M$ is compact.  Babillot
\cite[Thm.~1]{Babillot02b} showed that if $m_F$ is finite, then it is
mixing for the geodesic flow on $T^1M$ if the length spectrum of $\Ga$
is nonarithmetic (that is, if the set of translation lengths of the
elements of $\Ga$ is not contained in a discrete subgroup of
$\RR$). This condition, conjecturally always true, is known, for
example, if $\Ga$ has a parabolic element, if $\Lambda\Ga$ is not
totally disconnected (hence if $M$ is compact), or if $\wt M$ is a
surface or a (rank-one) symmetric space, see for instance
\cite{Dalbo99,Dalbo00}.

For every subset $A$ of $\wt M$ and every point $x$ in $\wt M
\cup \partial_\infty \wt M$, the {\em shadow of $A$ seen from $x$} is
the set $\OOO_xA$ of points at infinity of the geodesic rays or lines
starting from $x$ and meeting $A$. By Mohsen's shadow lemma (see
\cite[Lem.~3.10]{PauPolSha}), for every $x\in \wt M$, if $R>0$ is
large enough, there exists $c=c(R)>0$ such that for every $\ga\in\Ga$,
we have
\begin{equation}\label{eq:mohsen}
\frac{1}{c}\;e^{\int_x^{\ga y}(\wt F-\delta_{\Ga,\,F})}
\leq \mu^F_{x}(\OOO_x(B(\ga y,R)))\leq
c\; e^{\int_x^{\ga y}(\wt F-\delta_{\Ga,\,F})}\;. 
\end{equation}

Here is a new consequence of Mohsen's shadow lemma which will be
useful in this paper. Recall that a discrete group $G$ of isometries
of $\wt M$ is {\it convex-cocompact} if its limit set $\Lambda G$
contains at least two points, and if the action of $G$ on the convex
hull $\C\Lambda G$ in $\wt M$ of $\Lambda G$ has compact quotient.

\blemm\label{lem:zeromeaslimitsetsubgroup} Let $\Ga_0$ be a
convex-cocompact subgroup of $\Ga$ such that $\delta_{\Ga_0,\,F_0}
<\delta_{\Ga,\,F}$, where $F_0:\Ga_0\bs T^1 \wt M\ra \RR$ is the map
induced by $\wt F$.  Then $\mu^F_{x_0}(\Lambda\Ga_0)=0$.  
\elemm

\dem Since $\Ga_0$ is convex-cocompact, if $R$ is big enough, for
every $n\in\NN$, we have 
$$
\Lambda\Ga_0\subset \bigcup_{\ga\in\Ga_0,\;
  d(x_0,\,\ga x_0)\geq n}\;\OOO_{x_0}(B(\ga x_0,R))\;.
$$
Hence, by Equation \eqref{eq:mohsen}, there exists $c>0$ such that for
every $n\in \NN$,
$$
\mu^F_{x_0}(\Lambda\Ga_0)\leq \sum_{\ga\in\Ga_0,\;
  d(x_0,\,\ga x_0)\geq n}\mu^F_{x_0}(\OOO_{x_0}(B(\ga x_0,R)))
\leq c\sum_{\ga\in\Ga_0,\;
  d(x_0,\,\ga x_0)\geq n}\;e^{\int_{x_0}^{\ga x_0}(\wt F-\delta_{\Ga,\,F})}\;.
$$
The Poincaré series $Q_{\Ga_0,\,F_0,\,x_0}(\delta_{\Ga,\,F})$ converges,
as $\delta_{\Ga_0,\,F_0} <\delta_{\Ga,\,F}$. Since the remainder of a
converging series tends to $0$, this proves the result.  
\cqfd

\medskip A {\it parabolic} subgroup of $\Ga$ is a maximal infinite
subgroup $\Ga_0$ of $\Ga$ whose limit set $\Lambda\Ga_0$ is a
singleton. It is {\it bounded} if $\Ga_0\bs(\Lambda\Ga-\Lambda\Ga_0)$
is compact. For every bounded parabolic subgroup $\Ga_0$, if
$\Lambda\Ga_0=\{\xi_0\}$, there exists (see for instance
\cite{Bowditch95}) a unique $\Ga$-equivariant family $(\H_{\alpha
  \xi_0})_{\alpha \in \Ga/\Ga_0}$ of maximal closed horoballs in $\wt
M$ with pairwise disjoint interiors and with $\H_{\xi_0}$
centred at $\xi_0$. The horoball $\H_{\xi_0}$ is {\it precisely
  invariant} under $\Ga$, that is, its stabiliser in $\Ga$ is $\Ga_0$
and the inclusion $\stackrel{\circ}{\H_{\xi_0}}\;\subset \wt M$
induces an injection $\Ga_0\bs\!\stackrel{\circ}{\H_{\xi_0}} \;\ra
\Ga\bs \wt M$. Note that if $\Ga_0$ is a parabolic subgroup of $\Ga$
and if $m_F$ is finite, then we also have $\mu^F_{x_0} (\Lambda\Ga_0)
=0$ (see \cite[Lem.~5.15]{PauPolSha}).

\section{Hausdorff dimension of Patterson measures of 
potentials} 
\label{sec:ledrappier} 

Let $\wt M$ be a complete simply connected Riemannian manifold with
pinched negative sectional curvature at most $-1$.  Let $\Ga$ be a
nonelementary discrete group of isometries of $\wt M$. Let $\wt F:
T^1\wt M\to\RR$ be a H\"older-continuous $\Ga$-invariant function.
Assume that $\delta=\delta_{\Ga,\,F}$ is finite.  Let $\wt m_F$ be the
Gibbs measure on $T^1\wt M$ associated to a pair of 
Patterson densities $\big((\mu^{F\circ\iota}_x) _{x\in\wt M},
(\mu^F_{x})_{x\in\wt M} \big)$ for $(\Ga,F\circ\iota)$ and
$(\Ga,F)$. We use the notation introduced in Section
\ref{sec:reminder}. 

We fix in this section a point $x$ in $\wt M$. We denote by
$d_x$ the Gromov-Bourdon {\it visual distance} on $\partial_\infty \wt
M$ seen from $x$, defined (see \cite{Bourdon95}) by
\begin{equation}\label{eq:defidistvis}
d_x(\xi,\eta)=
\lim_{t\ra+\infty} e^{\frac{1}{2}(d(\xi_t,\,\eta_t)-d(x,\,\xi_t)-d(x,\,\eta_t))}\;,
\end{equation}
where $t\mapsto \xi_t,\eta_t$ are any geodesic rays ending at
$\xi,\eta$ respectively. We endow from now on $\partial_\infty \wt M$
with the distance $d_x$.

The aim of this section is to compute the Hausdorff dimension of the
Patterson measure $\mu^F_{x}$ associated to the potential $F$ (which
will be independent of $x$). Recall that the {\it Hausdorff dimension}
$\dim_H (\nu)$ of a finite nonzero measure $\nu$ on a locally compact
metric space $X$ is the greatest lower bound of the Hausdorff dimensions
$\dim_H (Y)$ of the Borel subsets $Y$ of $X$ with $\nu(Y) > 0$.

Let us give a motivation for such a computation. As mentioned in the
introduction, we are interested in this paper in studying whether the
set $E(\psi)$ of vectors of $T^1M$ that are well-spiraling, as
quantified by $\psi$, around a given closed geodesic $D_0$ has full or
zero measure for the Gibbs measure $m_F$. Varying the potential $F$
may be useful to estimate the Hausdorff dimension of $E(\psi)$: if
$\int_0^{+\infty} e^{\psi(t)(P(F_{|T^1D_0})-P(F))} dt$ diverges, as we
will prove in Section \ref{sec:spiral}, the set $E(\psi)$ has full
measure for $m_F$, and hence $\dim_H (E(\psi))\geq \dim_H(m_F)$. Note
that the Hopf parametrisation of $T^1\wt M$ is Hölder-continuous
(though usually not Lipschitz, except in particular when $\wt M$ is a
symmetric space), and $\wt m_F$ is in the same measure class as the
product measure $d\mu^{F\circ\iota}_x \otimes d\mu^F_{x}\otimes
dt$. Hence $\dim_H(m_F)$ may be estimated using $\dim_H(\mu^F_x)$
(that we will prove to be equal to $\dim_H(\mu^{F\circ\iota}_x)$), and
is in fact equal to $2\dim_H(\mu^F_x)+1$ if $\wt M$ is a symmetric
space.

The main result of this section, proving Theorem
\ref{theo:ledrappintro} in the introduction, is the following one. To
simplify the notation, let $h(m)=h_{\frac{m}{\|m\|}}(g^1)$ be the
(metric) entropy of the geodesic flow with respect to a finite nonzero
$(g^t)_{t\in\RR}$-invariant measure $m$ on $T^1M$ normalised to be a
probability measure.

\btheo\label{theo:ledrappier} If the Gibbs measure $m_F$ is finite and
if $F$ is $m_F$-integrable, then the Hausdorff dimension of the
Patterson measure $\mu^F_x$ on $(\partial_\infty \wt M,d_x)$
associated to $F$ satisfies
\begin{equation}\label{eq:computhausdimPatterson}
\dim_H(\mu^F_x)=\dim_H(\mu^{F\circ\iota}_x) = h(m_{F})\leq \delta_\Ga\;.
\end{equation}
If $M$ is convex-cocompact, then the last
inequality is an equality if and only if $F-P(F)$ is cohomologous to
the zero potential.  
\etheo

We think that the convex-cocompact assumption in the last claim could
be improved (see the comment at the end of this section).

The first claim is a generalisation of a result of Ledrappier
\cite{Ledrappier13}, who proved the theorem in the particular case
$F=0$.  Then $\mu_x^0$ is the standard Patterson measure of $\Ga$ and
the associated Gibbs measure $m_{F}$ is the Bowen-Margulis measure
$m_{\rm BM}$.  Let $\Lambda_c \Gamma$ denote the {\it conical (or
  radial) limit set}, that is, the set of $\xi \in \partial_\infty \wt
M$ for which $\liminf_{t \to +\infty} d(\rho(t), \Gamma x) < +\infty$,
where $\rho$ is any geodesic ray with point at infinity $\xi$. Let
$h_{\rm top} (g^1)$ be the topological entropy of the geodesic flow on
$T^1M$.  If $m_{\rm BM}$ is finite, then Ledrappier
\cite[Theo.~4.3]{Ledrappier13} proves furthermore that
$$
\dim_H(\mu^0_x)=h(m_{\rm BM})=h_{\rm top}
(g^1)=\dim_H(\Lambda_c\Ga)=\delta_\Ga\;.
$$
The second equality is due to Otal-Peigné \cite{OtaPei04}. The last
equality, which does not require the assumption that $m_{\rm BM}$ is
finite, is due to Bishop-Jones in constant curvature, to Hamenstädt
and to the first author (see \cite{Paulin97d})in general.

\medskip \dem Up to normalising $\mu^F_x$, which does not change its
Hausdorff dimension nor $\frac{m_F}{\|m_F\|}$, we may assume that
$\mu^F_x$ is a probability measure. The proof will follow from a
series of lemmas and propositions.  The following is a well known
useful alternative characterisation of the dimension of the measure,
which was also used by Ledrappier \cite[Prop.~2.5]{Ledrappier13}.

\blemm\label{lem:defdim} For any finite nonzero measure $\nu$ on a
compact metric space $X$, the Hausdorff dimension $\dim_H (\nu)$ is
the $\nu$-essential greatest lower bound on $x\in X$ of
$$
\liminf_{\epsilon\to 0} \;\frac{\ln \nu(B(x, \epsilon))}{\log \epsilon} \;.
$$
\elemm

For every $\xi \in \partial_\infty \widetilde M$, let $\rho_{\xi}:
[0,+\infty[\; \to \wt M$ be the geodesic ray with $\rho_\xi(0) = x$
and $\rho_\xi(+\infty) = \xi$.  The next lemma compares shadows of
balls in $\wt M$ with (visual) balls in $\partial_\infty \wt M$.

\blemm [Bourdon \cite{Bourdon95}]\label{lem:bourdon} For sufficiently
large $R>0$, there exists $D = D(R)$ such that, for all $\epsilon > 0$
and $\xi\in\partial_\infty\wt M$,
$$
\OOO_x(B(\rho_\xi(\log (1/\epsilon) + D), R)) 
\subset B_{d_x}(\xi, \epsilon ) 
\subset \OOO_x(B(\rho_\xi(\log (1/\epsilon) - D), R))\;.
$$
\elemm

Our first step in proving the theorem  is the following result.
 
\bprop\label{prop:majoparpression} 
\begin{enumerate}
\item[(1)] If $(\Ga, F)$ is of divergence type then $\dim_H (\mu_x^F)
  \leq \dim_H(\Lambda_c\Gamma) = \delta_\Gamma$;
\item[(2)] If $m_F$ is finite and if $F$ is $m_F$-integrable, then
  $\dim_H (\mu^F_x) \leq P(F) - \frac{1}{\|m_F\|}\,\int_{T^1M} F\,
  dm_F$.
\end{enumerate} 
\eprop

\dem By \cite[Theo.~5.12]{PauPolSha}, if $(\Ga, F)$ is of divergence
type, then the set $\Lambda_c\Gamma$ has full $\mu^F_x$-measure, and
thus the inequality in Part (1) follows immediately from the
definition of the Hausdorff dimension of measures.  The equality in
Part (1) has already been mentioned.

In order to prove Part (2), note that $(\Ga, F)$ is of divergence type
if $m_F$ is finite, by \cite[Coro.~5.15]{PauPolSha}. It hence suffices
by Lemma \ref{lem:defdim} to show that for $\mu^F_x$-almost every
$\xi$ in the full $\mu^F_x$-measure subset $\Lambda_c\Ga$, we have
$$
\liminf_{\epsilon \to 0} \;
\frac{\log \mu^F_x(B_{d_x}(\xi, \epsilon))}{\log \epsilon} \;\leq\; P(F) -
 \frac{1}{\|m_F\|}\, \int_{T^1M} F\, dm_F \;.
$$

Let $\xi\in\Lambda_c\Ga$. By the definition of $\Lambda_c\Ga$, there
exist $K\geq 0$, a sequence $(\ga_n)_{n\in\NN}$ in $\Ga$ and a
sequence $(t_n)_{n\in\NN}$ converging to $+\infty$ in $[0,+\infty[$
such that $d(\rho_\xi(t_n), \gamma_nx) \leq K$. By the triangle
inequality, we have $d(x,\ga_n x)\leq t_n+K$ and the ball
$B(\rho_\xi(t_n), R)$ contains the ball $B(\ga_nx, R - K)$, for every $R
\geq K$.  Let us apply the inclusion on the left in Lemma
\ref{lem:bourdon} with $\epsilon_n = e^{-t_n + D({R})}$, which tends
to $0$ as $n \to +\infty$ (hence in particular may be assumed to be in
$]0,1]$). We have
\begin{equation}\label{eq:majoepsilonn}
\frac{\log \mu^F_x(B_{d_x}(\xi, \epsilon_n))}{\log \epsilon_n}
\leq \frac{\log \mu^F_x(\OOO_x(B(\rho_\xi(t_n), R)))}{\log \epsilon_n}
\leq 
\frac{\log \mu^F_x(\OOO_x(B(\gamma_nx, R-K)))}{\log \epsilon_n}\;.
\end{equation}
By Mohsen's shadow lemma (see Equation \eqref{eq:mohsen}) and by
\cite[Theo.~6.1]{PauPolSha} which says that $P(F)=\delta_{\Ga,\,F}$, if $R$
is large enough, there exists $c>0$ such that, for every $n\in\NN$,
$$
\mu^F_{x}(\OOO_x(B(\ga_n x,R-K)))\geq 
\frac{1}{c}\;e^{\int_x^{\ga_n x}(\wt F-P(F))}\;.
$$
By Equation \eqref{eq:holdercontpotential}, we have $\int_x^{\ga_nx}
\wt F = \int_0^{t_n} \wt F(\dot\rho_\xi(s)) \,ds +
\operatorname{O}(1)$ as $n\ra+\infty$.  Thus Equation
\eqref{eq:majoepsilonn} gives, as $n\ra+\infty$,
\begin{align*}
\frac{\log \mu^F_x(B_{d_x}(\xi, \epsilon_n))}{\log \epsilon_n}&
\leq \frac{\int_x^{\gamma_n x} (\wt F - P(F))\; -\log c}{-t_n +D(R)} \\
&\leq \big(P(F) - \frac{1}{t_n} 
\int_0^{t_n} \wt F(\dot\rho_\xi(s)) \,ds\big)(1 + \operatorname{o}(1))\;. 
\end{align*}
By \cite[Theo.~5.4]{PauPolSha}, since $(\Ga, F)$ is of divergence
type, the geodesic flow in $T^1M$ is ergodic for $m_F$. Since $F$ is
$m_F$-integrable on $T^1M$, since $m_F$ is finite and by the
quasi-product structure of $\wt m_F$ in Hopf's parametrisation, for
$\mu^F_x$-almost every $\xi$, we have by Birkhoff's ergodic theorem
$$
\lim_{n\ra+\infty}\;\frac{1}{t_n} \int_0^{t_n} \wt F(g^s\dot\rho_\xi(0))\, ds 
\;=\; \frac{1}{\|m_F\|}\,\int_{T^1M} F\, dm_F\;.
$$
This proves Proposition \ref{prop:majoparpression}. \cqfd

\medskip
We next want to show that the reverse inequality holds. 
 
\bprop \label{prop:minoparpression} 
If $m_F$ is a finite measure and if $F$ is $m_F$-integrable, then
$\dim_H (\mu_x^F) \geq P(F) - \frac{1}{\|m_F\|}\,\int_{T^1M} F\, dm_F$.
\eprop

\dem To prove the result, by Proposition \ref{lem:defdim}, we only
need to show that for $\mu_x^F$-almost every $\xi$, we have 
$$
\liminf_{\epsilon \to 0}\; 
\frac{\log \mu^F_x(B_{d_x}(\xi, \epsilon))}{\log \epsilon} 
\;\geq\; P(F) -\frac{1}{\|m_F\|}\, \int_{T^1M} F\, dm_F\;.
$$

As in the proof of \cite[Prop.~4.6]{Ledrappier13}, since $m_F$ is
finite and by the quasi-product structure of $\wt m_F$, by
Poincar\'e's recurrence theorem and Birkhoff's ergodic theorem, for
$\mu^F_x$-almost every $\xi$, there exist $K>0$, a sequence
$(\ga_n)_{n\in\NN}$ in $\Ga$ and an increasing  sequence $(t_n)_{n\in\NN}$,
converging to $+\infty$ in $[0,+\infty[$, such that $d(\rho_\xi(t_n),
\gamma_nx) \leq K$, and such that the limit $\lim_{n\ra+\infty} t_n/n$
exists and is positive.

Let $R$ be big enough and let $c=c(R+K)$ be as in Mohsen's shadow lemma
(see Equation \eqref{eq:mohsen}), so that, for every $n\in\NN$,
$$
\mu^F_x(\OOO_x (B(\ga_n x, R + K)))\leq 
c\; e^{\int_x^{\gamma_nx} (F- P(F))}\;.
$$
By the triangle inequality, the ball $B(\ga_nx, R + K)$ contains the
ball $B(\rho_\xi(t_n), R)$. For every $n\in\NN$, let $\epsilon_n =
e^{-t_n - D(R)}$, which decreases to $0$. For every $\epsilon\in \;
]0,1]$ small enough, let $n=n(\epsilon)\in\NN$ be such that
$\epsilon_n \geq \epsilon > \epsilon_{n+1}$. By the inclusion on the
right in Lemma \ref{lem:bourdon} and by the same arguments as in the
end of the proof of the previous proposition, we have
\begin{align*}
\frac{\ln\mu^F_x(B_{d_x}(\xi, \epsilon))}{\ln \epsilon} & \geq
\frac{\ln \mu^F_x(B_{d_x}(\xi, \epsilon_n))}{\ln \epsilon_{n+1}} \geq 
\frac{\ln \mu^F_x(\OOO_x (B(\rho_\xi(t_n), R)))}{\ln \epsilon_{n+1}}
\\ &=\frac{-\ln \mu^F_x(\OOO_x (B(\rho_\xi(t_n), R)))}{t_{n+1} + D(R)}
 \geq \frac{-\ln \mu^F_x(\OOO_x (B(\ga_nx, R + K)))}{t_{n+1} + D(R)} 
  \\ &\geq \frac{-\int_x^{\ga_n x} (F- P(F))\;-\ln c}{t_{n+1} + D(R)} 
\geq  \frac{t_nP(F) - 
\frac{t_n+\operatorname{o}(t_n)}{\|m_F\|}\int_{T^1M} F \;dm_F 
+\operatorname{O}(1)}{t_{n+1} + D(R)} \;. 
 \end{align*}
 Taking the inferior limit as $\epsilon\ra 0$, since
 $\lim_{n\ra+\infty}\frac{t_n}{t_{n+1}}=1$, the result follows. \cqfd

\medskip Now, by the Variational Principle
\cite[Theo.~6.1]{PauPolSha}, since $m_F$ is finite and since $F$ is
$m_F$-integrable, we have $P(F)= h(m_F)+ \frac{1}{\|m_F\|}\,
\int_{T^1M} F\, dm_F$. Since $\iota:T^1M\ra T^1M$ conjugates
$(g^t)_{t\in\RR}$ to $(g^{-t})_{t\in\RR}$, and since
$m_{F\circ\iota}=\iota_*m_F$, we have $h(m_{F\circ\iota})=h(m_F)$.
Hence Equation \eqref{eq:computhausdimPatterson} in Theorem
\ref{theo:ledrappier} follows from Propositions
\ref{prop:majoparpression} and \ref{prop:minoparpression} applied to
both $F$ and $F\circ\iota$.

If $\Ga$ is convex-cocompact, then $m_F$ and $m_{\rm BM}=m_0$ are
finite and $F$ is integrable for $m_F$ and $m_0$. By the uniqueness in
the Variational Principle (see \cite[Theo.~6.1]{PauPolSha}), if
$h(m_F)= \delta_\Ga=h(m_0)$, then $\frac{m_F}{\|m_F\|}=
\frac{m_0}{\|m_0\|}$. By the Hamenstädt-Ledrappier correspondence (see
\cite{Ledrappier95, Hamenstadt97, Schapira04a} and the following
proposition) saying that if $\Ga$ is convex-cocompact, the cohomology
class of a potential with zero pressure is determined by
its associated Gibbs measure, the last claim of Theorem
\ref{theo:ledrappier} follows. \cqfd

\medskip We end this section by a comment on the correspondence
between the potentials and their associated Patterson measures,
which will be used at the end of this paper.

\bprop[Hamenstädt-Ledrappier]\label{prop:hamled} If $\Ga$ is
convex-cocompact, the map $\wt F\mapsto \mu=\mu^F_{x_0}$ induces a
bijection from the set of $\Ga$-invariant Hölder maps $\wt
F:\wt\Omega\Ga\ra\RR$ with zero pressure $P(F)=0$, up to cohomologous
maps, to the set of measure classes of Hölder quasi-invariant measures
$\mu$ on $(\Lambda\Ga,d_x)$ endowed with the action of $\Ga$.
Furthermore, for every hyperbolic element $\ga\in\Ga$ with attractive
fixed point $\ga^+\in\Lambda\Ga$, the period of $\ga$ for $\wt F$
satisfies
\begin{equation}\label{eq:periodcocycle}
\operatorname{Per}_F(\ga)=\ln\frac{d(\ga^{-1})_*\mu}{d\,\mu}(\ga^+)\;.
\end{equation}  
\eprop 

\dem The reader who is not interested in seeing how this result can be
deduced from \cite{Ledrappier95} (whose arguments extend from the
cocompact to the convex-cocompact case, as observed in
\cite{Schapira04a}) may skip this proof.

Recall that $\partial_\infty\wt M$ has a unique Hölder structure such
that for every $x\in\wt M$, the map $v\mapsto v_+$ from $T^1_x\wt M$
to $\partial_\infty\wt M$ (whose inverse will be denoted by
$\xi\mapsto v_{x,\,\xi}$) is a Hölder homeomorphism.

The following definitions are are taken from \cite{Ledrappier95}. A
{\it Hölder cocycle} for the action of $\Ga$ on $\partial_\infty\wt M$
is a map $c:\Ga \times \Lambda\Ga\ra\RR$, which is Hölder-continuous
in the second variable, such that $c(\ga\ga',\xi)=c(\ga,\ga'\xi)+
c(\ga',\xi)$ for all $\ga, \ga' \in\Ga$ and $\xi\in\Lambda\Ga$. The
{\it period} for $c$ of a hyperbolic element $\ga$ of $\Ga$ is
$c(\ga,\ga^+)$, where $\ga^+$ is the attractive fixed point of
$\ga$. Two Hölder cocycles $c$ and $c'$ are {\it cohomologous} if
there exists a Hölder-continuous map $U: \Lambda\Ga\ra\RR$ such that
$c(\ga,\xi)-c'(\ga,\xi)= U(\ga\xi)- U(\xi)$ for all $\ga\in\Ga$ and
$\xi\in\Lambda\Ga$.  Given a Hölder quasi-invariant measure $\mu$, its
{\it associated} Hölder cocycle is $c_\mu:(\ga,\xi)\mapsto -\ln
\frac{d(\ga^{-1})_*\mu} {d\,\mu} (\xi)$. The verification that this is
indeed a Hölder cocycle is immediate.

\medskip Fix $x_0\in\wt\Omega\Ga$. Given a potential (that is, a
$\Ga$-invariant Hölder map) $\wt F:\wt\Omega\Ga\ra\RR$, the map
$c_F:(\ga,\xi)\mapsto C^F_\xi(\ga^{-1}x_0,x_0)$ is a Hölder cocycle
(see \cite[Prop.~3.5 (2)]{PauPolSha} for its Hölder-continuity, $\wt
F$ being bounded since $\Ga\bs \wt\Omega\Ga$ is compact). Hence, by
the definition of a Patterson density, given a potential $\wt
F:\wt\Omega\Ga\ra\RR$, the measure $\mu^F_{x_0}$ is a Hölder
quasi-invariant measure, whose associated Hölder cocycle is $c_F$.  If
two potentials $\wt F$ and $\wt F^*$ are cohomologous, then their
associated Hölder cocycles $c_F$ and $c_{F^*}$ are cohomologous: it is
easy to check that if $\wt G:\wt\Omega\Ga\ra\RR$ is H\"older-continuous,
$\Ga$-invariant, differentiable along every flow line, and satisfies
Equation \eqref{eq:cohomologouspotential}, then the map
$U:\Lambda\Ga\ra\RR$ defined by $\xi\mapsto \wt G(v_{x_0,\,\xi})$ is
Hölder-continuous and satisfies $c_{F^*}(\ga,\xi)-c_{F}(\ga,\xi)=
U(\ga\xi)- U(\xi)$ for all $\ga\in\Ga$ and $\xi\in\Lambda\Ga$.

Let us relate the periods of a potential $F$ to the periods of the
Hölder cocycle $c_F$. Let $\ga$ be a hyperbolic element of $\Ga$, with
translation axis $A_\ga$, translation length $ \ell(\ga)$ and
attractive fixed point $\ga_+$. By the $\Ga$-invariance and the
cocycle property of $C^F$, if $p$ is the closest point to $x_0$ on
$A_\ga$, we have $C^F_{\ga^+}(\ga^{-1}x_0,x_0)
=C^F_{\ga^+}(\ga^{-1}p,p)$. Hence, by the definition of $C^F$, with
$t\mapsto \xi_t$ the geodesic ray from $p$ to $\ga^+$, we have (note
that there are sign differences with \cite{Ledrappier95})
\begin{align}
c_F(\ga,\ga^+)&=C^F_{\ga^+}(\ga^{-1}x_0,x_0) =C^F_{\ga^+}(\ga^{-1}p,p)
=\lim_{t\ra+\infty}\int_{p}^{\xi_t} (\wt F-P(F))-
\int_{\ga^{-1} p}^{t} (\wt F-P(F))\nonumber\\ &
= -\int_{\ga^{-1} p}^{p} (\wt F-P(F))
=P(F)\,\ell(\ga)- \operatorname{Per}_F(\ga)\;.
\label{eq:periodperiod}
\end{align}

By \cite[Théo.~1.c]{Ledrappier95}, two Hölder quasi-invariant measures
have the same measure class if and only if their associated Hölder
cocycles are cohomologous, and this holds if and only if the periods
of these Hölder cocycles are the same. By Liv\v{s}ic's theorem
(see \cite[Rem.~3.1]{PauPolSha}, two potentials $\wt F$ and $\wt F^*$
are cohomologous if and only they have the same periods. By Equation
\eqref{eq:periodperiod}, the periods of two potentials $\wt F$ and $\wt
F^*$ with zero pressure are the same if and only if the periods of the
associated Hölder cocycles $c_F$ and $c_{F^*}$ are the same. Hence the
map which associates to the cohomology class of a potential $\wt F$
the measure class of the Hölder quasi-invariant measure $\mu^F_{x_0}$
is well-defined, and is injective. To prove that it is surjective, we
start with a Hölder quasi-invariant measure $\mu$, we consider its
associated Hölder cocycle $c_\mu$, the proof of
\cite[Théo.~3]{Ledrappier95} shows that there exists a potential $\wt
F$ such that the Hölder cocycle $c_F$ is cohomologous to $c_\mu$, and
we apply again \cite[Théo.~1.c]{Ledrappier95} to get that
$\mu^F_{x_0}$ and $\mu$ have the same measure class.  

In order to prove \eqref{eq:periodcocycle}, if $\wt F$ is a potential
with $P(F)=0$, we have, by Equation \eqref{eq:periodperiod},
$$
\ln\frac{d(\ga^{-1})_*\mu^F_{x_0}}{d\,\mu^F_{x_0}}(\ga^+)=-c_F(\ga,\ga^+)
=\operatorname{Per}_F(\ga)\;.\;\;\;\Box
$$

\medskip It would be interesting to know if one could remove the
assumption that $\Ga$ is convex-cocompact, up to adding the
requirements on $\wt F$ that $\delta_{\Ga,\,F}$ is finite and
$(\Ga,F)$ is of divergence type, and on $\mu$ that $\mu$ is
ergodic. This would improve correspondingly the last claim of Theorem
\ref{theo:ledrappier} and simplify the statement of the requirement on
the class of measures under consideration in Theorem
\ref{theo:diophantine}.

\section{Almost sure spiraling for Gibbs states} 
\label{sec:spiral}

We will study in this section the generic asymptotic penetration
properties of the geodesic lines, in a negatively curved simply
connected manifold, under a discrete group of isometries, of a tubular
neighbourhood of a convex subset with cocompact stabiliser, not only
as in \cite{HerPau10} for the Bowen-Margulis measure, but for any
Gibbs measure.

\medskip Let $(\wt M,\Ga,\wt F,(\mu^{F\circ\iota}_x) _{x\in\wt M},
(\mu^F_{x})_{x\in\wt M}, \wt m_F )$ be as in the beginning of
Section \ref{sec:ledrappier}, with $\delta=\delta_{\Ga,\,F}$ finite.
We again use the notation introduced in Section \ref{sec:reminder}.

Recall that a subgroup $H$ of a group $G$ is {\it almost malnormal}
if, for every $g$ in $G-H$, the subgroup $gHg^{-1}\cap H$ is finite.
Let $\Ga_0$ be an almost malnormal and convex-cocompact subgroup of
$\Ga$, of infinite index in $\Ga$, let $C_0=\C\Lambda\Ga_0$ be the
convex hull of the limit set of $\Ga_0$. For instance, $C_0$ could be
the translation axis of a loxodromic element of $\Ga$, and $\Ga_0$ the
stabiliser of $C_0$ in $\Ga$ (see \cite[\S 4]{HerPau10} for an
explanation and for more examples). Up to adding assumptions on the
behaviour of the potential and on growth properties in cusp
neighbourhoods (including a gap property for the pressures), our
result should extend when $\Ga_0$ is assumed to be only geometrically
finite instead of convex-cocompact, or when $\Ga_0$ is a bounded
parabolic group (in which case $\Ga_0$ is also malnormal with infinite
index in $\Ga$) and $C_0$ is a precisely invariant closed horoball
centred at the singleton $\Lambda\Ga_0$. We restrict to the above case
for simplicity.

Let $F_0:\Ga_0\bs T^1\wt M\ra \RR$ be the map induced by $\wt F$, and
let $\delta_0=\delta_{\Ga_0,\,F_0}$ be the critical exponent of
$(\Ga_0,F_0)$. Note that $-\infty<\delta_0\leq \delta<+\infty$ by
\cite[Lem.~3.3 (iii)]{HerPau10}.

Let $\psi:[0,+\infty[\;\ra[0,+\infty[$ be a measurable map, such that
there exist $c_2,c_3>0$ such that for every $s,t\geq c_2$, if $s\leq
t+c_2$, then $\psi(s)\leq \psi(t)+c_3$.  Recall (see for instance
\cite[\S 5]{HerPau04}) that this condition is for instance satisfied
if $\psi$ is Hölder-continuous; it implies that $e^{\psi}$ is locally
bounded, hence it is locally integrable; and for every $\alpha>0$, the
series $\sum_{n\in\NN} e^{\alpha\,\psi(n)}$ converges if and only if
the integral $\int_0^{+\infty} e^{\alpha\,\psi(t)}\;dt$ converges.
Note that the constant $c_2$ and $c_3$ are unchanged by replacing
$\psi$ by $\psi+c$ for any $c\in\RR$.

Fix $\epsilon_0>0$. With the terminology of \cite{HerPau10}, let $\wt
E(\psi)$ be the set of {\it $(\epsilon_0,\psi)$-Liouville vectors for
  $(\Ga,\Ga_0)$} in $T^1\wt M$, that is, the set of $v\in T^1\wt M$
such that there exist sequences $(t_n)_{n\in\NN}$ in $[0,+\infty[$
converging to $+\infty$ and $(\ga_n)_{n\in\NN}$ in $\Ga$ such that for
every $t\in [t_n, t_n + \psi(t_n)]$, we have
$g_v(t)\in\ga_n\N_{\epsilon_0}C_0$. Note that $\wt E(\psi)$ is
invariant under the geodesic flow and under $\Ga$.

If $E$ is a set and $f,g:E\ra \;]0,+\infty[$ are maps, we write $f\,
\asymp\,g$ if there exists $c>0$ such that $\frac{1}{c} \,f\leq g\leq
c\,f$. The aim of this section is to prove the following result.

\btheo\label{theo:maintechnicalspiralGibbs} Assume that the
measure $m_F$ is finite, and that there exists $\kappa>0$ such that
$\sum_{\ga\in\Ga\;:\;t\leq d(x,\ga y)< t+\kappa} \; e^{\int_{x}^{\ga
    y} \wt F} \; \asymp\;e^{t\;\delta}$ and
$\sum_{\alpha\in\Ga_0\;:\;t\leq d(x,\alpha y)< t+\kappa} \;
e^{\int_{x}^{\alpha y} \wt F}\; \asymp\; e^{t\;\delta_0}$.  If
$\int_0^{+\infty} e^{\psi(t)(\delta_0-\delta)} \; dt$ diverges
(resp.~converges) then $\wt m_F$-almost every (resp.~no) point of
$T^1\wt M$ belongs to $\wt E(\psi)$.  
\etheo

\noindent {\bf Remarks. } (1) If the length spectrum of $\Ga$ is
nonarithmetic, then as said in Section \ref{sec:reminder}, the measure
$m_F$ is mixing for the geodesic flow on $T^1M$, hence by
\cite[Coro.~9.7]{PauPolSha}, we have $\sum_{\ga\in\Ga\;:\;d(x,\,\ga y)
  \leq t} \; e^{\int_{x}^{\ga y} \wt F}\; \sim c\;e^{t\;\delta}$ as
$t\ra + \infty$, for some $c>0$, a stronger requirement than the first
asymptotic hypothesis. Similarly, if the length spectrum of $\Ga_0$ is
nonarithmetic (this implies that $\Ga_0$ is nonelementary), then the
Gibbs measure $m_{F_0}$ of $(\Ga_0,F_0)$, being finite since $\Ga_0$
is convex-cocompact, is mixing, and the second asymptotic hypothesis
holds. The fact that the second asymptotic hypothesis holds when
$\Ga_0$ is elementary (that is, when $C_0$ is the translation axis of
a loxodromic element of $\Ga$) is given by \cite[Lem.~3.3
(ix)]{PauPolSha}.
 
\medskip (2) The above theorem implies Theorem \ref{theo:mainintro} in
the introduction. Indeed, $M$ being compact, the measure $m_F$ is
finite and the length spectrum of $\Ga$ is nonarithmetic. Hence the
two asymptotic hypotheses of Theorem
\ref{theo:maintechnicalspiralGibbs} hold by the previous remark. Note
that if $C_0$ is the translation axis of a loxodromic element of
$\Ga$, if $D_0$ is its image by $\wt M\ra M$, then $\delta_0=P(F_{\mid
  T^1D_0})$ by \cite[Lem.~3.3 (ix)]{PauPolSha}. We have $\delta=P(F)$
by \cite[Theo.~6.1]{PauPolSha}. Hence the conclusion of Theorem
\ref{theo:maintechnicalspiralGibbs} does imply Theorem
\ref{theo:mainintro}.

\bigskip\noindent{\bf Proof of Theorem
  \ref{theo:maintechnicalspiralGibbs}. } Before starting this proof,
let us give more informations on $\Ga_0$. Recall that $C_0$ is a
non-compact, closed convex subset of $\wt M$ such that:
\begin{enumerate}
\item[(1)] $C_0$ is $\Ga_0$-invariant and $\Ga_0\backslash C_0$ is
  compact; up to replacing $\Ga_0$ by $\operatorname{Stab}_{\Ga} C_0$,
  in which $\Ga_0$ has finite index and which remains almost malnormal
  (see the caracterisation \cite[Prop.~2.6 (3)]{HerPau10}), so that
  $\delta_0$ and the validity of the second asymptotic hypothesis of
  Theorem \ref{theo:maintechnicalspiralGibbs} are unchanged, we may
  and we will assume that $\Ga_0= \operatorname{Stab}_{\Ga} C_0$;
\item[(2)] by \cite[Prop.~2.6 (2), (4)]{HerPau10}, the limit set
  $\Lambda\Ga_0$ is {\it precisely invariant} (that is, $\ga
  \Lambda\Ga_0\cap\Lambda\Ga_0=\emptyset$ for every $\ga\in\Ga-\Ga_0$),
  and there exists $\kappa_0>0$ such that for every $\ga\in\Ga-\Ga_0$,
  the diameter of $\N_{\epsilon_0}C_0\cap \ga \N_{\epsilon_0}C_0$ is
  at most $\kappa_0$.
\end{enumerate}

\blemm\label{lem:gap}
If $\Ga_0$ is a convex-cocompact subgroup of $\Ga$, then $\delta_0< \delta$.
\elemm

\dem Since the Gibbs measure $m_{F_0}$ is finite, by
\cite[Coro.~6.1]{PauPolSha}, the probability measure
$\frac{m_{F_0}}{\|m_{F_0}\|}$ is an equilibrium state for the
potential $F_0$ on $\Ga_0\bs T^1\wt M$, whose support is contained in
the nonwandering set $\Omega\Ga_0=\Ga_0\bs \wt \Omega\Ga_0$ of the
geodesic flow on $\Ga_0\bs T^1\wt M$. Since $\Ga_0$ is malnormal in
$\Ga$, the canonical map $p:\Ga_0\bs T^1\wt M\ra \Ga\bs T^1\wt M$,
when restricted on the nonwandering sets, is a finite-to-one map, by
the above property (2). Hence $p_*\big(\frac{m_{F_0}}{\|m_{F_0}\|}
\big)$ is another equilibrium state for $F$ on $\Ga\bs T^1\wt M$. But
by \cite[Coro.~6.1]{PauPolSha}, this equilibrium state is unique,
hence $p_*\big(\frac{m_{F_0}}{\|m_{F_0}\|}\big)=\frac{m_{F}}{\|m_{F}\|}$. 

Since $\Ga_0$ is convex-cocompact and has infinite index in $\Ga$, its
limit set $\Lambda\Ga_0$ is a precisely invariant (by the above
property (2)) nonempty closed subset with empty interior in
$\Lambda\Ga$. Hence $\Ga\Lambda\Ga_0$ is a proper subset of
$\Lambda\Ga$ by Baire's theorem. Therefore the support of
$p_*m_{F_0}$, which is the image by $p$ of $\Ga_0\bs \wt\Omega\Ga_0
=\Ga_0\bs\{v\in T^1\wt M\;:\; v_-,v_+\in\Lambda\Ga_0\}$, is a proper
subset of the support $\Ga\bs \wt\Omega\Ga =\Ga\bs\{v\in T^1\wt M\;:\;
v_-,v_+\in\Lambda\Ga\}$ of $m_{F}$, a contradiction. 
\cqfd

\medskip We start the proof of Theorem
\ref{theo:maintechnicalspiralGibbs} by two reductions of the
statement.

\medskip
(i) Up to adding a big enough constant to $\wt F$, which does not
change $\wt m_F$, nor $\delta_0-\delta$, nor the asymptotics of the
series in the above statement, we assume that $\delta_0>0$.  In
particular, $\delta$ is finite and positive.

(ii) Let $x_0\in C_0$ be a basepoint. Let $R_0 >0$ and let $\wt{U_0}=
\pi^{-1}(B(x_0, R_0))$ be the set of the unit tangent vectors in
$T^1\wt M$ based at a point at distance less that $R_0$ of $x_0$. If
$R_0$ is big enough, then $m_F(\wt{U_0})>0$. Since $m_F$ is finite, it
is ergodic under the action of the geodesic flow on $T^1 M$ (see
\cite[Coro.~5.15]{PauPolSha}). Hence the result is equivalent to
proving that, when $R_0$ is big enough, if $\int_0^{+\infty}
e^{\psi(t)(\delta_0-\delta)} \; dt$ diverges (resp.~converges) then
$\wt m_F$-almost every (resp.~no) point of $\wt{U_0}$ belongs to $\wt
E(\psi)\cap\wt{U_0} $.

\medskip We now define the various subsets of $\wt{U_0}$ that we will
study during the proof of Theorem \ref{theo:maintechnicalspiralGibbs}.

Let $E_0$ be the set of $[\ga]\in \Ga/\Ga_0$ such that $d(x_0,\ga
C_0)\leq R_0+\epsilon_0$. Since $\Ga$ is discrete, and since $\Ga_0$
acts cocompactly on $C_0$, only finitely many distinct images of $C_0$
under $\Ga$ meet a given compact subset of $\wt M$. In particular, the
set $E_0$ is finite.

Since $\Ga_0\backslash C_0$ is compact, let $\Delta_0>0$ be such that
the restriction to the ball $B(x_0,\Delta_0)$ of the canonical
projection $C_0\ra\Ga_0\backslash C_0$ is onto. Choose and fix once and
for all a representative $\ga$ of $[\ga]\in \Ga/\Ga_0-E_0$ such that
if $p_\ga$ is the closest point to $x_0$ on $\ga C_0$, then
$d(p_\ga,\ga x_0)\leq \Delta_0$. We will use this representative
whenever a coset is considered. For every $[\ga]\in \Ga/\Ga_0-E_0$,
define
$$
D_\gamma=d(x_0,\ga C_0)=d(x_0,p_\ga)>0\;.
$$

\brema\label{rem:lem41HP} 
{\rm Note that by an argument similar to \cite[Lem.~4.1]{HerPau10}, for
  every $\lambda\in\RR$, there are only finitely many $[\ga]\in
  \Ga/\Ga_0-E_0$ such that $D_\gamma\leq\lambda$.}  
\erema

\blemm\label{lem:minocroisscoset} Assume that there exists $\kappa>0$
such that
$$
\sum_{\ga\in\Ga\;:\; t\leq d(x_0,\,\ga x_0)< t+\kappa} \; e^{\int_{x_0}^{\ga x_0} \wt F}\; 
\asymp\; e^{\delta\;t}\;\;\;{\rm and} \;\;\;
\sum_{\alpha\in\Ga_0\;:\;t\leq d(x_0,\,\alpha x_0) < t+\kappa} \;
e^{\int_{x_0}^{\alpha x_0} \wt F}\; \asymp\; e^{\delta_0\;t}\;.
$$ 
Then there exists $\kappa'\geq 1$ such that
$\sum_{[\ga]\in\Ga/\Ga_0\;:\; t\leq D_\ga < t+\kappa'} \;
e^{\int_{x_0}^{\ga x_0} \wt F}\; \asymp\; e^{\delta\;t}$.
\elemm

\dem We start by proving that there exist $c_4,c_5>0$ such
that for every $([\ga],\alpha)\in\Ga/\Ga_0\times \Ga_0$, we have
\begin{equation}\label{eq:controlheight}
D_\ga\leq d(x_0,\ga x_0)\leq D_\ga+\Delta_0\;,
\end{equation}
\begin{equation}\label{eq:distanceadditivity}
d(x_0,\ga x_0)+d(x_0,\alpha x_0)-c_4\leq d(x_0,\ga\alpha x_0)
\leq d(x_0,\ga x_0)+d(x_0,\alpha x_0)\;,
\end{equation}
\begin{equation}\label{eq:potentieladditivity}
\Big|\int_{x_0}^{\ga\alpha x_0} \wt F\;-
\int_{x_0}^{\ga x_0} \wt F\; -\int_{x_0}^{\alpha x_0} \wt F\Big|\leq c_5\;.
\end{equation}

Equation \eqref{eq:controlheight}, as well as the inequality on the
right hand side of Equation \eqref{eq:distanceadditivity}, follow by the
triangle inequality:
$$
D_\ga=d(x_0,\ga C_0)\leq d(x_0,\ga x_0)\leq d(x_0, p_\ga)+ 
d(p_\ga, \ga x_0)\leq D_\ga+\Delta_0\;.
$$

By the convexity of $\ga C_0$, the angle at $p_\ga$ of the geodesic
segments $[p_\ga,x_0]$ and $[p_\ga, \ga\alpha x_0]$ (if they are
non-trivial) is at least $\frac{\pi}{2}$. By hyperbolicity, the point
$p_\ga$ is hence at distance at most $\log(1+\sqrt{2})$ from a point
in $[x_0, \ga\alpha x_0]$. Thus $\ga x_0$ is at distance at most
$\Delta_0+ \log(1+\sqrt{2})$ from a point $u$ in $[x_0, \ga\alpha
x_0]$. By the triangle inequality, the inequality on the left hand
side of Equation \eqref{eq:distanceadditivity} follows with
$c_4=2(\Delta_0+\log(1+\sqrt{2}))$.

Let us apply Equation \eqref{eq:holdercontpotential} twice, with
$x=u$, $y=\ga x_0$ and with either $z=x_0$ or $z=\ga\alpha x_0$. Since
$d(\ga x_0,u)\leq \Delta_0+ \log(1+\sqrt{2})$, Equation
\eqref{eq:potentieladditivity} follows with
$$c_5=
2\big(c_1\,e^{\Delta_0+\ln(1+\sqrt{2})} +(\Delta_0+\ln(1+\sqrt{2}))
\max_{\pi^{-1}(B(x_0,\,\Delta_0+\ln(1+\sqrt{2})))} |\wt F|\,\big)\;.
$$

We are now going to use the following lemma.

\blemm \cite[Lem.~3.3]{HerPau10} \label{lem:herpautroitroi} For all
$A,\delta_0,\delta>0$, there exists $N\in\NN$ and $B>0$ such that for
all sequences $(a_k)_{k\in\NN}$ and $(b_k)_{k\in\NN}$ such that
$a_n\leq A\,e^{\delta\,n}$,  $b_n\leq A\,e^{\delta_0\,n}$ and
$\sum_{k=0}^n a_kb_{n-k}\geq \frac{1}{A}\,e^{\delta\,n}$ for every
$n\in\NN$ big enough, we have $\sum_{k=0}^N a_{n+k}\geq
B\,e^{\delta\,n}$ for every $n\in\NN$. 
\elemm

By the first asymptotic assumption in Lemma \ref{lem:minocroisscoset},
there exists $c>0$ such that, for every $t\geq \kappa$,
$$ 
\sum_{\ga\in \Ga,\;t-\kappa\leq d(x_0,\,\ga x_0)<t} 
e^{\int_{x_0}^{\ga x_0}\wt F}\geq \frac{1}{c}\;e^{\delta \;t}\;.
$$
We will use Lemma \ref{lem:herpautroitroi} by taking, for every
$k\in\NN$,
$$
a_k=\sum_{[\ga]\in \Ga/\Ga_0,\;k\leq D_\ga <k+1} e^{\int_{x_0}^{\ga
    x_0}\wt F}
\;\;\;{\rm and} \;\;\;
b_k=\sum_{\alpha\in \Ga_0,\;k-\kappa-c_4\leq d(x_0,\,\alpha x_0) <k}
e^{\int_{x_0}^{\alpha x_0}\wt F}\;.
$$
By Equation \eqref{eq:controlheight} and by the first asymptotic
assumption in Lemma \ref{lem:minocroisscoset}, there exists $C'>0$
such that, for every $k\in\NN$,
\begin{equation}\label{eq:controlupbk}
a_k\leq \sum_{\ga\in \Ga,\;k\leq d(x_0,\,\ga x_0)<k+1+\Delta_0} 
e^{\int_{x_0}^{\ga x_0}\wt F}\leq c'\;e^{\delta \;k}\;.
\end{equation}
By the second asymptotic assumption in Lemma
\ref{lem:minocroisscoset}, there exists $c''>0$ such that, for every
$k\in\NN$,
$$
b_k\leq c''\,e^{\delta_0 \,k}\;.
$$

Let $n\geq \kappa+c_4$ and $([\ga],\,\alpha)\in\Ga/\Ga_0\times \Ga_0$
satisfy $n-\kappa-c_4\leq d(x_0,\ga' x_0)<n-c_4$ where $\ga'= \ga
\alpha$. Let $k=\lfloor D_\ga\rfloor$ be the integral part of $D_\ga$.
By Equation \eqref{eq:distanceadditivity}, we hence have
$$
0\leq k\leq D_\ga\leq d(x_0,\ga x_0)\leq d(x_0,\ga\alpha x_0)- 
d(x_0,\alpha x_0)+c_4\leq n\;,
$$
and
\begin{multline*}
n-\kappa-c_4-k\leq d(x_0,\ga\alpha x_0)- d(x_0,\ga x_0)\\\leq  
d(x_0,\alpha x_0)\\ \leq
d(x_0,\ga\alpha x_0)- d(x_0,\ga x_0)+c_4\leq n-k\;.
\end{multline*}
Therefore, respectively by the definition of $a_k$ and $b_{n-k}$, and
by Equation \eqref{eq:potentieladditivity}, we have
\begin{align*}
\sum_{k=0}^n a_kb_{n-k} &
=\sum_{k=0}^n\;\;\sum_{[\ga]\in \Ga/\Ga_0,\;k\leq D_\ga
  <k+1} \;\;\sum_{\alpha\in \Ga_0,\;n-k-\kappa-c_4\leq d(x_0,\,\alpha x_0) <n-k}
e^{\int_{x_0}^{\ga x_0}\wt F}e^{\int_{x_0}^{\alpha x_0}\wt F}
\\ &
\geq e^{-c_5}\sum_{\ga'\in \Ga,\;n-\kappa-c_4\leq d(x_0,\,\ga' x_0)<n-c_4} 
e^{\int_{x_0}^{\ga' x_0}\wt F}\geq e^{-c_5}\;\frac{1}{c}\; e^{\delta \;(n-c)}\;.
\end{align*}
Applying Lemma \ref{lem:herpautroitroi} with $A=\max\{c',c'',
c\,e^{c_5+\delta c_4}\}$ gives the lower bound required to prove Lemma
\ref{lem:minocroisscoset}.  The upper bound follows from Equation
\eqref{eq:controlupbk}. \cqfd

\bigskip
For every $r>0$ and $\beta\in\Ga$, let
$$
A_{\beta}(r) = \{v\in \wt{U_0}\;:\; g_v([0,+\infty[)\cap B(\beta x_0,r)\neq
\emptyset\}\;.
$$
Let us fix a positive constant $c_6\geq \kappa$ (depending only on
$\epsilon_0,\Delta_0, R_0, \kappa$ and $\psi$) to be made precise
later on. For every $k\in\NN$, define $I_k$ to be the set of $[\ga]\in
\Ga/\Ga_0$ such that $k\leq D_\ga<k+1$, and let $J_k=J_k(\psi)$ be the
set of pairs $([\ga],\alpha)\in \Ga/\Ga_0\times \Ga_0$ such that
$k\leq D_\ga<k+\kappa'$ (where $\kappa'$ is given by Lemma
\ref{lem:minocroisscoset}) and $\psi(k)\leq d(x_0,\alpha x_0)<
\psi(k)+c_6$. For every $k\in \NN$, let
$$
A_k(r,\psi)= \bigcup_{([\ga],\,\alpha)\in J_k}
\;\;A_{\ga\alpha}(r)\;.
$$

These sets are related to the set $\wt E(\psi)$ that we want to study
by the following result.
Recall that if $(B_k)_{k\in\NN}$ is a sequence of subsets of a given
set, one defines $\limsup_k B_k= \bigcap_{n\in\NN}\bigcup_{k\geq n}
B_k$.

\bprop\label{prop:geometrictranslation} If $r\geq \epsilon_0+
\Delta_0$, there exist $c'_5,c''_5>0$ such that, up to sets of $\wt
m_F$-measure zero,
$$
\limsup_{k} A_k(r,\psi+c''_5)\subset \wt E(\psi)\cap\wt{U_0} \;,
$$
and if $\psi(t)\geq c'_5$ for $t$ big enough,
$$
\wt E(\psi)\cap(\wt{U_0}-Z) \subset \limsup_{k} A_k(r,\psi-c'_5)\;.
$$
\eprop

\dem Let us first prove the second inclusion. Let $c_0=
\epsilon_0+2\operatorname{arsinh}(\coth \epsilon_0)$. Let $c'_0=
c_3\big\lceil \frac{2R_0+c_0}{c_2}\big\rceil$, with $c_2,c_3$ the
constants appearing in the assumption on $\psi$. Let $c'_5=
\epsilon_0+2\Delta_0 +R_0+c_0+c'_0$. Assume that $\psi(t)\geq c'_5$
for $t$ big enough.

Let $v\in \wt E(\psi)\cap\wt{U_0}$. For every $n\in\NN$, there exist
sequences $(t_n)_{n\in\NN}$ in $[0,+\infty[$ converging to $+\infty$
and $([\ga_n])_{n\in\NN}$ in $\Ga/\Ga_0$ such that for every $t\in
[t_n, t_n + \psi(t_n)]$, we have $g_v(t)\in
\ga_n\N_{\epsilon_0}C_0$. Let $n\in \NN$. The geodesic line $g_v$
enters in $\ga_n\N_{\epsilon_0}C_0$ at a time $t_n^-$ at most $t_n$.
Up to extracting a subsequence, we may assume, by Remark
\ref{rem:lem41HP}, that $[\ga_n]\notin E_0$, so that
$D_{\ga_n}>R_0+\epsilon_0$ and $t_n^->0$.

Let $k_n=\lfloor D_{\ga_n}\rfloor$. Let us prove that $k_n\ra+\infty$
as $n\ra+\infty$, up to sets of $\wt m_F$-measure zero of elements $v
\in\wt E(\psi)\cap\wt{U_0}$. Otherwise, up to extracting a
subsequence, $(\ga_n)_{n\in\NN}$ is constant by Remark
\ref{rem:lem41HP}.  Hence $v_+$ belongs to the set
$\ga_0 \partial_\infty C_0$ of accumulation points of $ \ga_0 C_0$ in
$\partial_\infty \wt M$. By Lemma \ref{lem:zeromeaslimitsetsubgroup},
the $\mu^F_{x_0}$-measure of $\partial_\infty C_0=\Lambda\Ga_0$ is
zero. Hence, since the action of $\Ga$ preserves the sets of
$\mu^F_{x_0}$-measure zero by the properties of the Patterson
densities, we have $\mu^F_{x_0}\Big(\bigcup_{\beta\in\Ga}
\beta \partial_\infty C_0\Big) =0$.  By the decomposition of $\wt m_F$
in Hopf's parametrisation (see Equation \eqref{eq:defigibbs}), the
$\wt m_F$-measure of the set of $v\in \wt E(\psi)$ such that $v_+\in
\bigcup_{\beta\in\Ga} \beta \partial_\infty C_0$ is zero. This proves
the above claim.

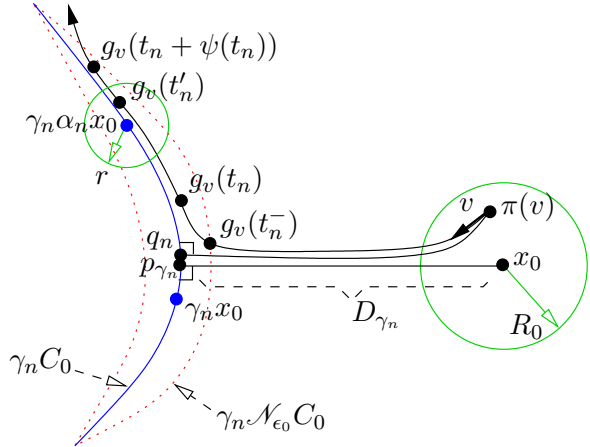
\begin{center}
\input{fig_rightinclusion.pstex_t}
\end{center}

\medskip Let $q_n$ be the closest point to $\pi(v)$ on $\ga_nC_0$. It
satisfies $d(p_{\ga_n},q_n)\leq R_0$, since closest point maps do not
increase the distances. Note that the point $q_n$ is at distance
$\epsilon_0$ from the entry point in $\ga_n\N_{\epsilon_0}C_0$ of the
geodesic segment from $\pi(v)$ to $q_n$. By the penetration properties
of geodesic rays in $\epsilon_0$-neighbourhoods of convex subsets of
$\operatorname{CAT}(-1)$ metric spaces (see
\cite[Lem.~2.3]{ParPau10GT}), we have $d(q_n, g_v(t^-_n)) \leq c_0=
\epsilon_0+2\operatorname{arsinh}(\coth \epsilon_0)$.  Hence, by the
triangle inequality,
\begin{equation}\label{eq:distpied}
d(\ga_n x_0, g_v(t^-_n))\leq 
d(\ga_n x_0, p_{\ga_n})+d(p_{\ga_n},q_n)+d(q_n,g_v(t^-_n))\leq 
\Delta_0+R_0+c_0\;.
\end{equation}
Again by the triangle inequality, we  have
\begin{equation}\label{eq:majokn}
k_n\leq D_{\ga_n}=d(x_0,p_{\ga_n})\leq 
d(x_0,\pi(v))+ t^-_n+d(g_v(t^-_n),q_n)+d(q_n,p_{\ga_n})\leq 
t_n +2R_0+c_0\;.
\end{equation}
Up to extracting a subsequence, we may assume that $\psi(k_n)\geq
c'_5$ and that $t_n,k_n\geq c_2$. By the assumption on $\psi$ and
since $c'_0= c_3\big\lceil \frac{2R_0+c_0}{c_2}\big\rceil$, we have by
Equation \eqref{eq:majokn},
$$
t^-_n+\psi(k_n)\leq t_n + \psi(t_n) + c'_0\;.
$$
Let $t'_n=t^-_n+\psi(k_n)-c'_0$, which belongs to $[t^-_n, t_n+
\psi(t_n)]$, since $\psi(k_n)\geq c'_5\geq c'_0$. By convexity, the
point $g_v(t'_n)$ belongs to $\N_{\epsilon_0}C_0$.  Let $\alpha_n$ be
an element of $\Ga_0$ such that
\begin{equation}\label{eq:controltprime}
  d(g_v(t'_n),\ga_n\alpha_n x_0)\leq \epsilon_0+\Delta_0\;,
\end{equation}
which exists by the definition of $\Delta_0$.  By the triangle
inequality, and by Equation \eqref{eq:distpied}, we have
\begin{align*}
\big|\;d(\ga_n x_0,\ga_n\alpha_n x_0)-
d(g_v(t^-_n),g_v(t'_n))\;\big|&\leq 
d(g_v(t'_n),\ga_n\alpha_n x_0) + 
d(\ga_n x_0,g_v(t^-_n))\\ & \leq \epsilon_0+2\Delta_0 +R_0+c_0\;.
\end{align*}
Hence
\begin{align}
\big|\;d(x_0,\alpha_n x_0)-\psi(k_n)\;\big|&
=|d(\ga_n x_0,\ga_n\alpha_n x_0)-|t'_n-t^-_n|-c'_0|
\nonumber\\ & \leq \epsilon_0+2\Delta_0 +R_0+c_0+c'_0=c'_5
\label{eq:alphaPsi}\;.
\end{align}
Define  $c_6=\max\{2c'_5, \kappa\}$ (which only depends on $\epsilon_0,
\Delta_0, R_0, \kappa$ and $\psi$). Assume that $r\geq \epsilon_0+
\Delta_0$.  For every $n \in \NN$, we hence have $v\in
A_{\ga_n\alpha_n}(r)$ by Equation \eqref{eq:controltprime}. Besides,
$([\ga_n],\alpha_n)\in J_{k_n}(\psi-c'_5)$ since $k_n=\lfloor
D_{\ga_n}\rfloor$ and $\kappa'\geq 1$, and by Equation
\eqref{eq:alphaPsi}.  Therefore $v\in A_{k_n}(r,\psi-c'_5)$. This
proves the second inclusion in Proposition
\ref{prop:geometrictranslation}.

\bigskip Let us now prove the first inclusion. By hyperbolicity and an
argument of (strict) convexity (see for instance
\cite[Lem.~2.2]{ParPau10GT}), there exists $c''_0=c''_0(\epsilon_0,r)$
such that if a geodesic segment has endpoints at distance at most
$\max\{R_0,r\} +\ln(1+\sqrt{2})$ from two points in $C_0$ at distance
at least $c''_0$ one from the other, then this geodesic segment enters
$\N_{\epsilon_0} C_0$. Let $c''_5=\max\big\{c''_0+\Delta_0,c_3+\Delta_0+R_0+2c_0+r+
\big\lceil\frac{2R_0+c_0+1}{c_2}\big\rceil\big\}$.

Let $v\in \wt U_0$ and let $(k_n)_{n\in\NN}$ be a sequence in $\NN$
converging to $+\infty$.  Assume that $v\in A_{k_n}(r,\psi+c''_5)$ for
every $n$ in $\NN$. Let $([\ga_n],\alpha_n)\in J_{k_n}(\psi+c''_5)$ be
such that $v\in A_{\ga_n\alpha_n}(r)$: there exists $\tau_n\geq 0$
such that $g_v(\tau_n)\in B(\ga_n\alpha_n x_0,r)$. Since $d(\pi(v),
x_0)\leq R_0$, by the properties of closest point projections in
CAT$(-1)$-space, there exists $\tau'_n\in [0,\tau]$ such that
$d(g_v(\tau'_n), p_{\ga_n})\leq \max\{R_0,r\} +\ln(1+\sqrt{2})$. By
the definition of $c''_0$ and since $d(p_{\ga_n}, \ga_n\alpha_n x_0)
\geq d(\ga_n x_0,\ga_n\alpha_n x_0)-\Delta_0\geq \psi(k_n)+ c''_5-
\Delta_0\geq c''_0$, the geodesic line $g_v$ enters $\ga_n
\N_{\epsilon_0} C_0$.

\begin{center}
\input{fig_leftinclusion.pstex_t}
\end{center}

Let $t_n^-$ be the entry time of $g_v$ inside $\ga_n \N_{\epsilon_0}
C_0$, which satisfies, by Equation \eqref{eq:distpied},
$$
d(g_v(t^-_n),\ga_n x_0)\leq \Delta_0+R_0+c_0\;.
$$
Let $t^+_n$ be either $\tau_n$ if $g_v(\tau_n)\in \ga_n
\N_{\epsilon_0} C_0$ or the exit time of $g_v$ out of $\N_{\epsilon_0}
C_0$ otherwise. Again by \cite[Lem.~2.3]{ParPau10GT} and since closest
point maps do not increase the distances, if $q'_n$ is the closest
point to $g_v(\tau_n)$ on $C_0$, we have
$$
d(g_v(t^+_n),\ga_n \alpha_nx_0)\leq 
d(g_v(t^+_n),q'_n)+d(q'_n,\ga_n \alpha_n x_0)\leq c_0+r\;.
$$
As in Equation \eqref{eq:majokn}, we have $k_n\geq D_{\ga_n}-1\geq
t^-_n -2R_0-c_0-1$. Hence
$$
\psi(k_n)\geq \psi(t^-_n) -
c_3\Big\lceil\frac{2R_0+c_0+1}{c_2}\Big\rceil\;.
$$
By the triangle inequality and since $([\ga_n],\alpha_n)\in J_{k_n}
(\psi+c''_5)$, we have
\begin{align*}
t^+_n-t^-_n&\geq d(\ga_n x_0,\ga_n\alpha_n x_0)
-d(g_v(t^-_n),\ga_n x_0)-d(g_v(t^+_n),\ga_n \alpha_nx_0)\\ &\geq
\psi(k_n) +c''_5-  \Delta_0-R_0-2c_0-r\\ &\geq
\psi(t^-_n)\;,
\end{align*}
where $c''_5$ is a constant depending only on $\Delta_0, R_0,
\epsilon_0,r$ and $\psi$. Hence $v$ belongs to $\wt E(\psi)$,
which proves the result. \cqfd

\bigskip 
In a series of lemmae and propositions, we now state the required
properties of the sets $A_{\ga\alpha}(r)$ for $([\ga],\alpha)\in
\Ga/\Ga_0\times \Ga_0$ and $A_k(r,\psi)$ for $k\in\NN$.

\medskip We start by the following estimate on the mass of the
$A_{\ga\alpha}(r)$'s. Before stating it, let us motivate it. Let $d'$
be the distance on $T^1\wt M$ induced by Sasaki's Riemannian metric on
$T\wt M$ (when $\Ga$ is cocompact, any Riemannian distance on $T^1\wt
M$ is allowed). Recall that, for $\epsilon>0$ and $T\geq 0$, the {\it
  dynamical $(\epsilon,T)$-ball} centred at a point $v\in T^1\wt M$
is
$$
B_{\epsilon,\,T}(v)=\{w\in T^1\wt M\;:\; \forall\;t\in[0,T],\;\;
d'(\flow t w,\flow t v)\leq \epsilon\}\;.
$$ 
We proved in \cite[Prop.~3.16]{PauPolSha} (which was in fact written
after the first version of this paper), using a minor modification of
these dynamical balls, that Gibbs measure satisfy the Gibbs property
(when $\Ga$ is torsion free and cocompact, see for instance
\cite[Theo.~3.3]{BowRue75} for the lower bound, and
\cite[Lem.~20.3.4]{KatHas95} in the discrete time case): for every
$\epsilon>0$, for all $v\in T^1\wt M$ and $T\geq 0$ such that $v,\flow
T(v)$ map to a given compact subset of $\Ga\backslash T^1\wt M$, we
have
$$
\wt m_F\big(B_{\epsilon,\,T}(v)\big) \asymp 
e^{\int_0^T \wt F(\flow tv)\;dt\;-T\;P(F)}\;.
$$ 
Now, $A_{\ga\alpha}(r)$ is almost such a dynamical ball. Indeed, let
$v_{\ga\alpha}$ be the unit tangent vector at $x_0$ of the geodesic
segment from $x_0$ to $\ga\alpha x_0$, and let $T_{\ga\alpha} = d(x_0,
\ga\alpha x_0)$ (see the figure below). Our set $A_{\ga\alpha} (r)$
contains $B_{\epsilon_-,T_{\ga\alpha}} (v_{\ga\alpha})$ and is
contained in $B_{\epsilon_+,T_{\ga\alpha}}(v_{\ga\alpha})$ for some
positive constants $\epsilon_\pm$ depending only on $R_0,r$. The
following result (or rather Equation \eqref{eq:dynamicalballvolume})
is hence closely related to this Gibbs property.

\bprop\label{prop:coroshadowlem} If $r$ and $R_0$ are big enough, there
exists $c_7=c_7(r)>0$ such that for all but finitely many
$([\ga],\alpha)\in \Ga/\Ga_0\times \Ga_0$, we have
$$
\frac{1}{c_7}\; e^{\int_{x_0}^{\ga x_0}(\wt F-\delta)}
\;e^{\int_{x_0}^{\alpha x_0} (\wt F-\delta)} \leq \wt m_F 
\big(A_{\ga\alpha}(r)\big) 
\leq c_7\; e^{\int_{x_0}^{\ga x_0}(\wt F-\delta)}
\;e^{\int_{x_0}^{\alpha x_0} (\wt F-\delta)}\;.
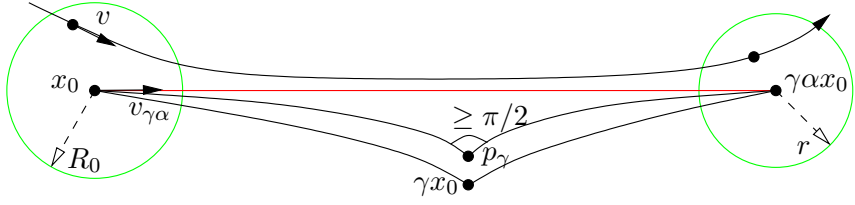
$$
\eprop

\begin{center}
\input{fig_dynamicalball.pstex_t}
\end{center}

\medskip 
\dem 
For every $R_0>0$ and $\beta\in\Ga$, define $B_{R_0,\,\beta}=
\bigcap_{z'\in B(\beta x_0,\, r)} \OOO_{z'} B(x_0,\frac{R_0}{3})$,
which is contained in (and is a pertubation of) the shadow
$\OOO_{\beta x_0} B(x_0,\frac{R_0}{3})$ (see the picture below). Since
$\Ga$ is nonelementary, the support of the Patterson measures is not
reduced to one point, hence $m= \inf_{\xi\in\partial_\infty \wt M}
\;\|\mu^{F\circ\iota}_{x_0}\|- \mu^{F\circ\iota}_{x_0} (\{\xi\})$ is
positive. By hyperbolicity (as first remarked by Sullivan), for every
$\xi\in\partial_\infty \wt M$, the family $\big(\,^c\OOO_\xi
B(x_0,R)\big)_{R>0}$ is a fundamental system of neighbourhoods of
$\xi$. By compacteness and discreteness, there exists hence $R_0>0$
such that for all but finitely many $\beta \in\Ga$, we have
\begin{equation}\label{eq:minoPattball}
\mu^{F\circ\iota}_{x_0}\big(B_{R_0,\,\beta}\big)\geq \frac{m}{2}\;.
\end{equation}
By the definition of $A_\beta(r)$, the set of points $v_+$ for $v$ in
$A_\beta(r)$ is exactly
$$
A_\beta(r)_+=\bigcup_{z\in B(x_0, \,R_0)} \OOO_z B(\beta x_0,r)\;,
$$
which is a bit larger than the shadow $\OOO_{x_0} B(\beta x_0,r)$ (see
the picture below).  By a minor modification of Mohsen's shadow lemma
(see Equation \eqref{eq:mohsen}), we have, if $r$ is big enough,
\begin{equation}\label{eq:mesPatAbetr}
\mu^F_{x_0}(A_\beta(r)_+)\asymp e^{\int_{x_0}^{\beta x_0} (\wt F
  -\delta)}\;.
\end{equation}

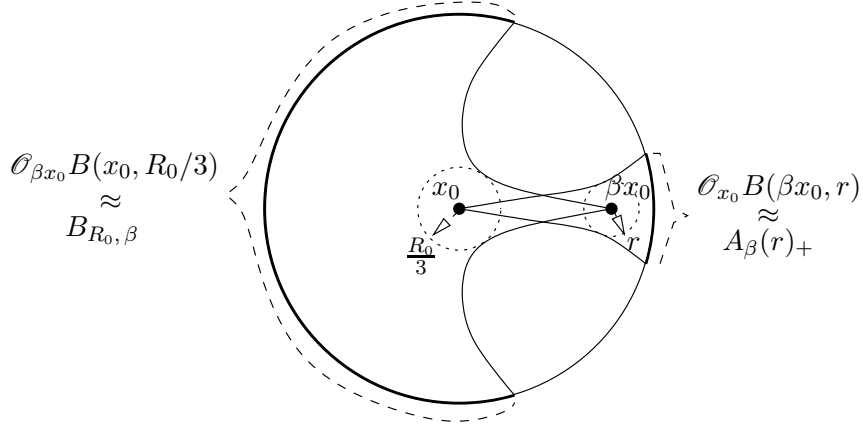
\begin{center}
\input{fig_bishadows.pstex_t}
\end{center}

\medskip 
Let us first prove, using Hopf's parametrisation defined by
the point ${x_0}$ and the definition of $A_\beta(r)$, that for all but
finitely many $\beta\in\Ga$, we have
\begin{equation}\label{eq:encadrementbinclusion}
B_{R_0,\,\beta} \times A_\beta(r/2)_+\times [-R_0/3,R_0/3] \subset
A_\beta(r)\subset
\partial_\infty\wt M\times A_\beta(r)_+\times [-R_0,R_0]\;.
\end{equation}
The inclusion on the right hand is immediate.  To prove the other one,
let $v\in T^1\wt M$ be such that if $p$ is the closest point to $x_0$
on the geodesic ray $g_v([0,+\infty[)$, then $v_-\in B_{R_0,\,\beta}$, $v_+\in
A_\beta(r/2)_+$ and $d(\pi(v),p)\leq R_0/3$ (see the picture
below). Let us prove that $v\in A_\beta(r)$.

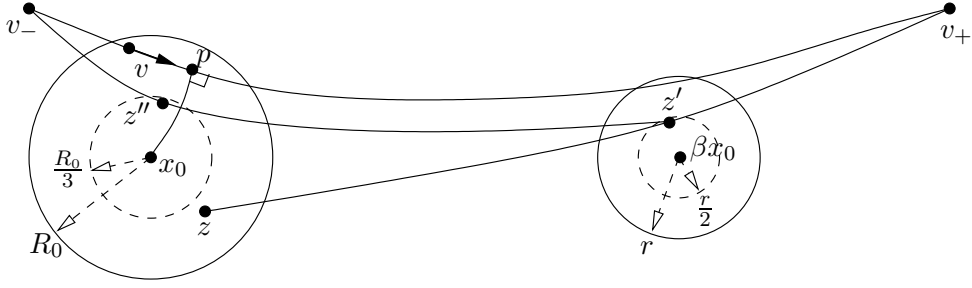
\begin{center}
\input{fig_minovoldynball.pstex_t}
\end{center}

By the definition of $A_\beta(r/2)_+$, let $z\in B(x_0,R_0)$ be such
that the geodesic ray $[z,v_+[$ meets $B(\beta x_0,r/2)$ at a point
$z'$. By the definition of $B_{R_0,\,\beta}$, the geodesic ray
$[z',v_-[$ meets $B(x_0,R_0/3)$ at a point $z''$.  Since $d(z,z'')\leq
2R_0$, and since by discreteness, for all but finitely many $\beta$ in
$\Ga$, the distance $d(z,z')$ is big, the angle at $z'$ between the
geodesic segments $[z',z]$ and $[z',z'']$ is small, hence the angle at
$z'$ between the geodesic rays $[z',v_-[$ and $[z',v_+[$ is close to
$\pi$. Hence the geodesic line $g_v$ between $v_-$ and $v_+$ is close
to the union of these two rays. In particular, since $g_v$ passes
close to $z'\in B(\beta x_0,r/2)$, it enters the ball $B(\ga x_0,r)$,
and since it passes close to $z''\in B(x_0,\frac{R_0}{3})$, the point
$p$ belongs to $B(x_0,\frac{R_0}{2})$ and hence $\pi(v)$ belongs to
$B(x_0,R_0)$, which proves the result.

\medskip Now, for every $v\in A_\beta(r)$, since $d(\pi(v),{x_0})\leq
R_0$, the point ${x_0}$ is at distance at most $R_0$ from a point on
the geodesic line between the endpoints $v_-, v_+$. Hence by Equation
\eqref{eq:holdercontcocycle}, there exists $c'_6\geq 0$ (depending
only on $R_0$, on $\max_{\wt{U_0}}|\wt F|<+\infty$, on the H\"older
constants of $\wt F$ and on the bounds of the sectional curvature of
$\wt M$) such that, for every $v\in A_\beta(r)$, we have
$$
-c'_6\leq C^{F\circ\iota}_{v_-}(x_0,\,\pi(v)),\;
C^F_{v_+}(x_0,\,\pi(v))\leq c'_6\;.
$$  
Therefore, by definition of the Gibbs measure $\wt m_F$, we have,
using Equations \eqref{eq:minoPattball} and
\eqref{eq:encadrementbinclusion},
$$
e^{-2\,c'_6}\;\frac{m}{2}\;\mu^F_{x_0}(A_\beta(r/2)_+)\;(2R_0/3)\leq 
\wt m_F\big(A_\beta(r)\big) \leq 
e^{2\,c'_6}\;\|\mu^{F\circ\iota}_{x_0}\|\;\mu^F_{x_0}( A_\beta(r)_+)\;(2R_0)\;.
$$
Hence, by Equation \eqref{eq:mesPatAbetr}, for some constant
$c''_6\geq 1$, we have
\begin{equation} \label{eq:dynamicalballvolume}
\frac{1}{c''_6}\;e^{\int_{x_0}^{\beta x_0} (\wt F -\delta)}\leq
\wt m_F\big(A_\beta(r)\big)\leq 
c''_6\;e^{\int_{x_0}^{\beta x_0} (\wt F -\delta)}\;.
\end{equation}

For all $[\ga]\in \Ga/\Ga_0$ and $\alpha\in \Ga_0$, by Equations
\eqref{eq:distanceadditivity} and \eqref{eq:potentieladditivity}, if
$\beta=\ga\alpha$, we have
$$
\Big|\int_{x_0}^{\beta x_0} (\wt F -\delta)\;-
\int_{x_0}^{\ga x_0} (\wt F -\delta)\;
-\int_{x_0}^{\alpha x_0} (\wt F -\delta)\Big|\leq c_5+\delta \,c_4\;.
$$
The result follows. \cqfd


\medskip
Note, as it will be important later on, that the contributions of
$\ga$ and of $\alpha$ are decoupled in this Proposition
\ref{prop:coroshadowlem}.

The sets $A_{\ga\alpha}(r)$ satisfy the following almost disjointness
property in shells.

\blemm\label{lem:almostdisjoint} For every $r>0$, there exists $c_8=
c_8(r)>0$ such that for every $k\in\NN$, for every subset $P$ of $J_k=
J_k(\psi)$,
$$
\frac{1}{c_8}\sum_{([\ga],\,\alpha)\in P} \wt m_F
\big(A_{\ga\alpha}(r)\big)\leq \wt m_F \big(\bigcup_{([\ga],\,\alpha)\in
  P} A_{\ga\alpha}(r)\big)\leq \sum_{([\ga],\,\alpha)\in P} \wt m_F
\big(A_{\ga\alpha}(r)\big)\;,
$$
and for every subset $Q$ of $I_k$
$$
\frac{1}{c_8}\sum_{[\ga]\in Q} \wt m_F
\big(A_{\ga}(r)\big)\leq \wt m_F \big(\bigcup_{[\ga]\in
  Q} A_{\ga}(r)\big)\leq \sum_{[\ga]\in Q} \wt m_F
\big(A_{\ga}(r)\big)\;.
$$
\elemm

\dem The inequality on the right hand side of the first claim is
immediate. In order to obtain the one on the left hand side, let us prove
that there exists $c_8\in\NN-\{0\}$ such that for all $k\in\NN$ and
$v\in T^1\wt M$, the number of $([\ga],\,\alpha)\in J_k$ such that
$v\in A_{\ga\alpha}(r)$ is at most $c_8$, which implies the result.

Let $([\ga],\,\alpha),([\ga'],\,\alpha')\in J_k$ such that $v\in
A_{\ga\alpha}(r)\cap A_{\ga'\alpha'}(r)$.  By Equation
\eqref{eq:controlheight} and \eqref{eq:distanceadditivity}, and by the
definition of $J_k$, we have
$$
k+\psi(k)-c_4\leq D_\ga+\psi(k)-c_4\leq d(x_0,\ga x_0)+d(x_0,\alpha
x_0)-c_4 \leq d(x_0,\ga\alpha x_0)
$$ and
$$
d(x_0,\ga\alpha x_0)\leq
 d(x_0,\ga x_0)+d(x_0,\alpha x_0) \leq D_\ga+\Delta_0+\psi(k)+c_6\leq
k+\kappa'+\Delta_0+\psi(k)+c_6\,.
$$
Similarly $k+\psi(k)-c_4\leq d(x_0,\ga'\alpha' x_0)\leq
k+\psi(k)+\kappa'+\Delta_0+c_6$.  

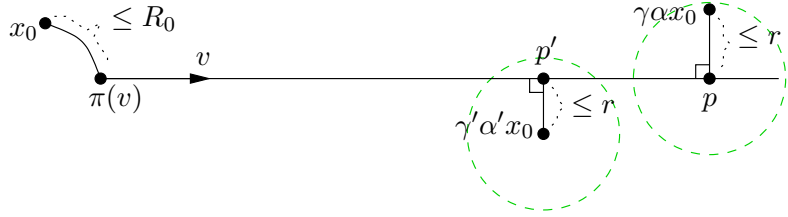
\begin{center}
\input{fig_brochette.pstex_t}
\end{center}

Let $p$ and $p'$ be the closest points on $g_v([0,+\infty[)$ to $\ga\alpha
x_0$ and $\ga'\alpha' x_0$ respectively. They satisfy $d(p,\ga\alpha
x_0), d(p',\ga'\alpha' x_0)\leq r$ since $v\in A_{\ga\alpha}(r)\cap
A_{\ga'\alpha'}(r)$. We may assume, up to permuting $\ga\alpha$ and
$\ga'\alpha'$, that $p'$ belongs to the geodesic segment
$[\pi(v),p]$. Since closest point maps do not increase the distances,
by the triangle inequality, and since $v\in \wt U_0$, we have
\begin{align*}
d(p,p') &=d(p,\pi(v))-d(\pi(v),p') \leq 
d(\ga\alpha x_0,\pi(v))-d(\pi(v),p')\\ &
\leq d(\ga\alpha x_0,x_0)+d(x_0,\pi(v)) - 
d(\ga'\alpha' x_0,x_0) +d(\pi(v),x_0)+d(p',\ga'\alpha' x_0)\\ &
\leq (k+\psi(k) +\kappa'+\Delta_0+c_6)+R_0-(k+\psi(k)-c_4)+R_0+r\;.
\end{align*}
Hence, again by the triangle inequality,
$$
d(\ga\alpha x_0,\ga'\alpha' x_0)\leq \kappa'+\Delta_0+c_6+ 2R_0 +c_4+3r\;.
$$
Now the first claim follows from the discreteness of $\Ga$, which
implies that there are only finitely many elements $\beta$ in $\Ga$
such that $\beta x_0$ belongs to a ball of centre $x_0$ with given
radius.

\medskip
The second claim is proven similarly.  \cqfd

\medskip The two results above allow to estimate the mass of the
$A_k(r,\psi)$'s, as follows.

\bprop\label{prop:estimatemeasure} Assume that there exists $\kappa>0$
such that
$$
\sum_{\ga\in\Ga\;:\; t\leq d(x,\ga y)< t+\kappa} 
\; e^{\int_{x}^{\ga y} \wt F}\; \asymp\;e^{\delta\;t}\;\;\;
{\rm and} \;\;\; \sum_{\alpha\in\Ga_0\;:\;t\leq d(x,\,\alpha y) 
< t+\kappa} \; e^{\int_{x}^{\alpha y} \wt F}\; \asymp\;
e^{\delta_0\;t}\;.
$$ 
If $r$ is big enough, there exists $c_9>0$ such that, for every $k\in
\NN$, we have
$$
\frac{1}{c_9}\; e^{\psi(k)(\delta_0-\delta)}\leq 
\wt m_F \big( A_k(r,\psi)\big)\leq 
c_9\; e^{\psi(k)(\delta_0-\delta)}\;.
$$
\eprop

It follows from this proposition and the assumption on the function
$\psi$ that the series $\sum_{k\in\NN}\wt m_F \big( A_k(r,\psi)\big)$
converges if and only if the integral $\int_0^{+\infty}
e^{\psi(t)(\delta_0-\delta)} \; dt$ converges.

\medskip \dem By Equation \eqref{eq:controlheight} and by the first
asymptotic assumption in the statement of Proposition
\ref{prop:estimatemeasure}, there exists $c>0$ such that for every
$k\in\NN$,
\begin{equation}\label{eq:controlupbkbis}
\sum_{[\ga]\in \Ga/\Ga_0,\;k\leq D_\ga<k+\kappa'} e^{\int_{x_0}^{\ga x_0}\wt F}
\leq \sum_{\ga\in \Ga,\;k\leq d(x_0,\,\ga x_0)<k+\kappa'+\Delta_0} 
e^{\int_{x_0}^{\ga x_0}\wt F}\leq c\;e^{\delta \;k}\;.
\end{equation}
By the second asymptotic assumption in Proposition
\ref{prop:estimatemeasure}, since $c_6\geq \kappa$, there exists
$c'>0$ such that, for every $t\in[0,+\infty[\,$,
\begin{equation}\label{eq:controlupck}
\frac{1}{c'}\;e^{\delta_0 \;t}\leq \sum_{\alpha\in \Ga_0,\;t\leq d(x_0,\,\ga x_0)<t+c_6} 
e^{\int_{x_0}^{\alpha x_0}\wt F} \leq c'\;e^{\delta_0 \;t}\;.
\end{equation}

Let us first prove the inequality in the right hand side in
Proposition \ref{prop:estimatemeasure}. Let $r$ be big enough and
$k\in \NN$. Respectively by Lemma \ref{lem:almostdisjoint} with
$P=J_k$ and the definition of $A_k(r,\psi)$, by Proposition
\ref{prop:coroshadowlem}, by the definition of $J_k$, by Equations
\eqref{eq:controlupbkbis} and \eqref{eq:controlupck}, we have
\begin{align*}
  m_F\big(A_k(r,\psi)\big)& \leq \sum_{([\ga],\,\alpha)\in J_k} \wt
  m_F \big(A_{\ga\alpha}(r)\big)\leq c_7\;\sum_{([\ga],\,\alpha)\in
    J_k} e^{\int_{x_0}^{\ga x_0}(\wt F-\delta)}
  \;e^{\int_{x_0}^{\alpha x_0} (\wt F-\delta)}\\ & =
  c_7\;\sum_{[\ga]\in \Ga/\Ga_0,\;k\leq D_\ga<k+\kappa'} e^{\int_{x_0}^{\ga
      x_0}(\wt F-\delta)}\; \sum_{\alpha\in \Ga_0,\;\psi(k)\leq
    d(x_0,\,\alpha x_0)<\psi(k)+c_6} \;e^{\int_{x_0}^{\alpha x_0} (\wt
    F-\delta)} \\ & \leq c_7\,(c\;e^{\delta \;k})\,e^{-\delta \;k}
  (c'\;e^{\delta_0 \;\psi(k)})\,e^{-\delta \;\psi(k)}=
c_7\,c\,c'\, e^{\psi(k)(\delta_0-\delta)}\;.
\end{align*}
This proves the inequality on the right hand side in Proposition
\ref{prop:estimatemeasure}. 

\medskip Let us now prove similarly the inequality on the left hand
side in Proposition \ref{prop:estimatemeasure}.  By Lemma
\ref{lem:minocroisscoset}, there exists $c''>0$ such that, for every
$t\in[0,+\infty[\,$,
\begin{equation}\label{eq:minocountdepthww}
\sum_{[\ga]\in\Ga/\Ga_0\;:\; t\leq D_\ga < t+\kappa'} \;
e^{\int_{x_0}^{\ga x_0} \wt F}\; \geq\; \frac{1}{c''}\;e^{\delta\;t}\;.
\end{equation}
Respectively by Lemma \ref{lem:almostdisjoint} with $P=J_k$ and the
definition of $A_k(r,\psi)$, by Proposition \ref{prop:coroshadowlem},
by the definition of $J_k$, by Equations \eqref{eq:minocountdepthww},
\eqref{eq:controlheight} and \eqref{eq:controlupck}, we have
\begin{align*}
  m_F\big(A_k(r,\psi)\big)& \geq
  \frac{1}{c_8}\sum_{([\ga],\,\alpha)\in J_k} \wt m_F
  \big(A_{\ga\alpha}(r)\big)\geq
  \frac{1}{c_7c_8}\;\sum_{([\ga],\,\alpha)\in J_k} e^{\int_{x_0}^{\ga
      x_0}(\wt F-\delta)} \;e^{\int_{x_0}^{\alpha x_0} (\wt
    F-\delta)}\\ & = \frac{1}{c_7c_8}\;\sum_{[\ga]\in
    \Ga/\Ga_0,\;k\leq D_\ga<k+\kappa'} e^{\int_{x_0}^{\ga x_0}(\wt
    F-\delta)} \sum_{\alpha\in \Ga_0,\;\psi(k)\leq d(x_0,\,\alpha
    x_0)<\psi(k)+c_6} \;e^{\int_{x_0}^{\alpha x_0} (\wt F-\delta)} \\
  & \geq \frac{1}{c_7c_8}\,(\frac{1}{c''}\;e^{\delta \;k})\,e^{-\delta
    \;(k+\kappa'+\Delta_0)} (\frac{1}{c'}\;e^{\delta_0 \;\psi(k)})\,e^{-\delta
    (\psi(k)+c_6)}\\
  &= \frac{1}{c_7\,c_8\,c'\,c''\,e^{\delta(\kappa'+\Delta_0+c_6)}}\;
  e^{\psi(k)(\delta_0-\delta)}\;.
\end{align*}
This proves Proposition \ref{prop:estimatemeasure}. \cqfd

\medskip The following result is a quasi-independence property of the
sets $A_k(r,\psi)$ for $k\in\NN$.

\bprop\label{prop:quasiindep} Under the hypotheses of Proposition
\ref{prop:estimatemeasure}, there exists a constant $c_{10}>0$ such
that  for every $k\neq k'$ in $\NN$, if $\psi\geq c_{10}$, we have
$$
\wt m_F \big( A_k(r,\psi)\cap A_{k'}(r,\psi)\big)\leq c_{10}\;
\wt m_F \big( A_k(r,\psi)\big) \;\wt m_F \big( A_{k'}(r,\psi)\big)\;.
$$
\eprop

\medskip \dem The proof has two parts, a geometric one and a
measure-theoretic one. We state the geometric one as a lemma.

\blemm \label{lem:geomquasiindep} There exist $c'_{10}>0$ and $r'>r$
such that for every $k<k'$ in $\NN$, for every $([\ga],\alpha)\in J_k$
and $([\ga'],\alpha')\in J_{k'}$, if $\psi\geq c'_{10}$ and if
$A_{\ga\alpha}(r)$ meets $A_{\ga'\alpha'}(r)$, then $A_{\ga'}(r)$ is
contained in $A_{\ga\alpha}(r')$.  \elemm

\dem Let $k<k'$ in $\NN$ and $([\ga],\alpha)\in J_k$ and $([\ga'],
\alpha')\in J_{k'}$. If $A_{\ga\alpha}(r)\cap A_{\ga'\alpha'}(r)$ is
non empty, there exists $v\in \wt {U_0}$ such that $g_v(\RR)$ meets
$B(\ga\alpha x_0,r)$ and $B(\ga'\alpha' x_0,r)$. Let $q,q'$ be the
closest points to $\pi(v)$ on the convex sets $\ga C_0,\ga' C_0$. Let
$p,p'$ be the closest points to $q,q'$ on the geodesic ray
$g_v([0,+\infty[)$. Let $x,x'$ be the closest points to $\ga\alpha
x_0$ and $\ga'\alpha' x_0$ on $g_v([0,+\infty[)$

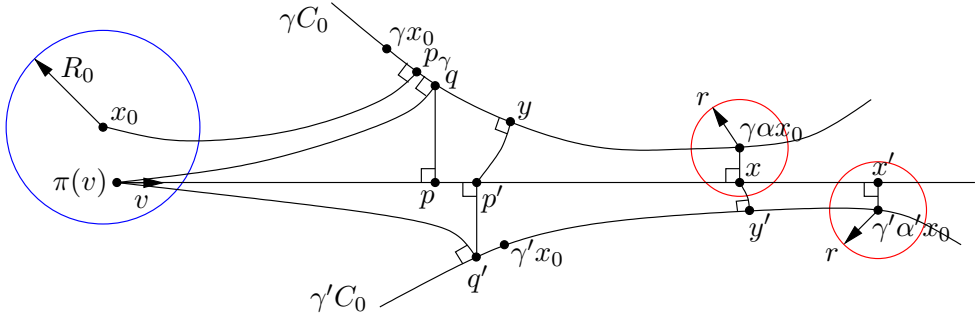
\begin{center}
\input{fig_quasiindependent.pstex_t}
\end{center}

We have $d(\ga\alpha x_0,x)\leq r$ and $d(\ga'\alpha' x_0,x')\leq r$.
By the properties of geodesic triangles in CAT$(-1)$-spaces and by the
convexity of $C_0$, we have $d(p,q),d(p',q')\leq r+\ln(1+\sqrt{2})$.
By the choice of the representatives of elements in $\Ga/\Ga_0$, the
closest point $p_\ga$ to $x_0$ on $\ga C_0$ is at distance at most
$\Delta_0$ from $\ga x_0$. Hence, since closest point maps do not
increase the distances and by the triangle inequality,
\begin{equation}\label{eq:majordpgax}
d(p,\ga x_0)\leq d(p,q)+d(q,p_\ga)+d(p_\ga, \ga x_0)\leq
r+\ln(1+\sqrt{2})+R_0+\Delta_0\;.
\end{equation}
Hence, by Equation \eqref{eq:controlheight} and the definition of
$J_k$, with $c=\kappa' +2R_0+2\Delta_0+r+\ln(1+\sqrt{2})$, we have
$$
d(\pi(v),p)\leq d(\pi(v),x_0)+d(x_0,\ga x_0)+d(\ga x_0,p)\leq 
R_0+(D_\ga+\Delta_0)+d(\ga x_0,p)\leq k+ c\;,
$$
and $d(\pi(v),p)\geq k-c$ by the inverse triangle inequality.
Similarly, 
$$
k'-c\leq d(\pi(v),p')\leq k'+ c\;.
$$
By similar arguments, if $c'= R_0+\Delta_0+2r+\ln(1+\sqrt{2})+c_6$, we
have
\begin{equation}\label{eq:majominodpxpsi}
\psi(k)-c'\leq d(p,x)\leq \psi(k)+ c'\;\;\;{\rm and}\;\;\;
\psi(k')-c'\leq d(p',x')\leq \psi(k')+ c'\;.
\end{equation}

Assume first that $\pi(v),x,p'$ are in this order on $g_v([0,
+\infty[)$.  Any geodesic ray, with origin at distance at most $R_0$
from $x_0$ and passing at distance at most $r$ from $\ga' x_0$, passes
at distance at most $2r+\ln(1+\sqrt{2})+R_0+\Delta_0$ from $p'$ by the
analog for $\ga'$ of Equation \eqref{eq:majordpgax}, hence by
convexity passes at distance at most $c''=\max\{2R_0,\;
2r+\ln(1+\sqrt{2})+R_0+\Delta_0\}$ from $x$, thus passes at distance
at most $c''+r$ from $\ga\alpha x_0$.  Therefore, if $r'\geq c''+r>r$,
then $A_{\ga'}(r)$ is contained in $A_{\ga\alpha}(r')$.

Assume now that $\pi(v),p',x$ are in this order on $g_v([0,+\infty[)$
(see the picture above). There exists a constant $c'_{10}>0$
(depending only the hyperbolicity constant $\ln(1+\sqrt{2}$ and on $r$)
such that if $\psi \geq c'_{10}$, then $\pi(v),p,x$ and $\pi(v),p',x'$
are in this order on $g_v([0,+\infty[)$.

Since $k'\geq k$, either $\pi(v),p,p'$ are in this order on $g_v([0,
+\infty[)$, or $p'\in[ \pi(v),p]$ is at distance at most $2c$ from
$p$, since then
$$
d(p,p')=d(p,x)-d(p',x)\leq (k+c)-(k'-c)\leq 2c\;.
$$
In both cases, by convexity, $p'$ is at distance at most $2c+r+ \ln(1+
\sqrt{2})$ from a point $y$ in $\ga C_0$.  Similarly, by Equation
\eqref{eq:majominodpxpsi} and since $\psi$ satisfies $\psi(s)\leq
\psi(t)+c_3$ if $s\leq t$, either $\pi(v),x,x'$ are in this order on
$g_v([0, +\infty[)$, or $x'\in [\pi(v),x]$ is at distance at most
$2c+2c'+c_3$ from $x$. In both cases, $x'$ is at distance at most
$2c+2c'+c_3+r+\ln(1+\sqrt{2})$ from a point $y'$ in $\ga' C_0$.

If for a contradiction $d(p',x)> R$ for arbitrarily large constants
$R$, then the geodesic segments $[y,\ga\alpha x_0]$ and $[q',y']$,
have endpoints at bounded distance from the long geodesic segment
$[p',x]$. Hence they have their endpoints at bounded distance while
being long, if $R$ is large. By hyperbolicity, this implies that
$\N_{\epsilon_0} (\ga C_0)\cap \N_{\epsilon_0} (\ga' C_0)$ contains a
long segment if $R$ is large. Taking $R$ large enough, this
contradicts the fact that the diameter of this intersection, since
$\ga\neq \ga$ in $\Ga/\Ga_0$, is at most the constant $\kappa_0$, as
explained in the beginning of the proof of Theorem
\ref{theo:maintechnicalspiralGibbs}.

Therefore $d(p',x)\leq R$ for some $R\geq 0$. Any geodesic ray, with
origin at distance at most $R_0$ from $x_0$ and passing at distance at
most $r$ from $\ga' x_0$, passes at distance from $\ga\alpha x_0$ at most
\begin{align*}
r+d(\ga'x_0,\ga\alpha x_0)&\leq r+d(\ga'x_0,p')+d(p',x)+d(x,\ga\alpha x_0)
\\ &\leq R+3r+\ln(1+\sqrt{2})+R_0+\Delta_0\;,
\end{align*}
by the analog for $\ga'$ of Equation \eqref{eq:majordpgax}.
Therefore, if $r'\geq R+3r+\ln(1+\sqrt{2})+R_0+\Delta_0>r$, then
$A_{\ga'}(r)$ is contained in $A_{\ga\alpha}(r')$.
\cqfd

\bigskip 
Now, let us use Lemma \ref{lem:geomquasiindep} to prove
Proposition \ref{prop:quasiindep}. Let $k,k'$ be elements of $\NN$
with $k<k'$. 

For every $([\ga],\alpha)\in J_k$, let $I_{[\ga],\,\alpha,\,k'}\subset
I_{k'}$ be the set of $[\ga']\in \Ga/\Ga_0$ such that there exists
$\alpha'\in\Ga_0$ with $([\ga'],\alpha')\in J_{k'}$ such that the
intersection $A_{\ga\alpha}(r)\cap A_{\ga'\alpha'}(r)$ is non empty.
Then respectively by Proposition \ref{prop:coroshadowlem}, by the
second part of Lemma \ref{lem:almostdisjoint} with $Q=
I_{[\ga],\,\alpha,\,k'}$, by Lemma \ref{lem:geomquasiindep} and the
definition of $I_{[\ga],\,\alpha,\,k'}$, and by Proposition
\ref{prop:coroshadowlem} (twice), we have
\begin{align}
\sum_{[\ga']\in I_{[\ga],\alpha,k'}} 
e^{\int_{x_0}^{\ga' x_0}(\wt F-\delta)}& \leq 
\sum_{[\ga']\in I_{[\ga],\alpha,k'}} 
c_7(r')\;\wt m_F \big(A_{\ga'}(r')\big)\nonumber \\ & \leq 
c_7(r')\,c_8(r')\;\wt m_F \Big(\bigcup_{[\ga']\in I_{[\ga],\alpha,k'}}
A_{\ga'}(r')\Big)\nonumber \\ & \leq 
c_7(r')\,c_8(r')\;\wt m_F \big(A_{\ga\alpha}(r')\big)
\nonumber \\ &
\leq c_7(r)\,c_7(r')^2\,c_8(r')\;\wt m_F \big(A_{\ga\alpha}(r)\big)
\label{eq:majotempotss}\;.
\end{align}

By the assumptions of Proposition \ref{prop:quasiindep}, there
exists $c>0$ such that for every $t\in\RR$,
\begin{equation}\label{eq:controlgammaO}
\sum_{\alpha'\in \Ga_0\;:\;d(x_0,\alpha' x_0)<t} 
\;e^{\int_{x_0}^{\alpha' x_0} \wt F}\;\leq c\;e^{t\;\delta_0}\;.
\end{equation}
To simplify the notation, let $A_k=A_k(r,\psi)$. Respectively by the
definition of $A_k$, by Proposition \ref{prop:coroshadowlem}, by Equation
\eqref{eq:controlgammaO}, by Equation \eqref{eq:majotempotss} with
$c'=c\;c_7(r)^2\,c_7(r')^2\,c_8(r')\;e^{c_6\delta_0}$, and by
Proposition \ref{prop:estimatemeasure} and Lemma
\ref{lem:almostdisjoint} with $P=J_k$, we have
\begin{align*}
  & \wt m_F(A_k\cap A_{k'})\\&\leq \sum_{([\ga],\,\alpha)\in J_k}\;\;
  \sum_{([\ga'],\,\alpha')\in J_{k'}\;:\;
A_{\ga\alpha}(r)\cap A_{\ga'\alpha'}(r)\neq \emptyset} 
\wt m_F\big(A_{\ga'\alpha'}(r)\big)\\
&\leq\sum_{([\ga],\,\alpha)\in J_k}\;\;\sum_{[\ga']\in I_{[\ga],\,\alpha,\,k'}}
\;c_7(r)\; e^{\int_{x_0}^{\ga' x_0}(\wt F-\delta)}
\sum_{\tiny\begin{array}{c}
\alpha'\in \Ga_0\\\psi(k')\leq d(x_0,\alpha' x_0)<\psi(k')+c_6
\end{array}} 
\;e^{\int_{x_0}^{\alpha' x_0} (\wt F-\delta)}\\ &\leq 
c\;c_7(r)\;e^{(\psi(k')+c_6)\delta_0-\delta\psi(k')}\;
\sum_{([\ga],\,\alpha)\in J_k}\;\sum_{[\ga']\in I_{[\ga],\,\alpha,\,k'}}\;
e^{\int_{x_0}^{\ga' x_0}(\wt F-\delta)}\\ &
\leq c'\;e^{\psi(k')(\delta_0-\delta)}\;
\sum_{([\ga],\,\alpha)\in J_k}\;\wt m_F \big(A_{\ga\alpha}(r)\big)
\\ &
\leq c'\; c_9\;c_8(r)\;\wt m_F (A_{k'})\;\wt m_F (A_{k})\;.
\end{align*}

This proves Proposition \ref{prop:quasiindep}.
\cqfd

\bigskip Let us now conclude the proof of Theorem
\ref{theo:maintechnicalspiralGibbs}.  The following version of the
Borel-Cantelli Lemma is well-known (see for instance
\cite{Sprindzuk69}).

\bprop \label{prop:borelcantelli}
Let $(Z,\nu)$ be a measured space with finite nonzero measure.  Let
$(A_n)_{n\in\NN}$ be a sequence of measurable subsets of $Z$ such that
there exists $c>0$ with $\nu(A_n\cap A_m)\leq c\;\nu (A_n) \,\nu(A_m)$
for all distinct $n,m$ in $\NN$. Then $\nu(\limsup_n A_n)>0$
if and only if the series $\sum_{n\in\NN}\nu(A_n)$ diverges. 
\eprop

We apply this result with $(Z,\nu)=(\wt{U_0},\wt m_{\mid \wt{U_0}})$,
which satisfies the hypothesis if $R_0$ is big enough as in the
reductions at the beginnning of the proof of Theorem
\ref{theo:maintechnicalspiralGibbs}. Let $r=\epsilon_0+ \Delta_0$ and
let $c'_5,c''_5$ be given by Proposition
\ref{prop:geometrictranslation}.

Assume first that the integral $\int_0^{+\infty} e^{\psi(t)(\delta_0-
  \delta)} \; dt$ diverges, which is still true if a constant is added
to $\psi$. The quasi-independence assumption of Proposition
\ref{prop:borelcantelli} is satisfied if $A_n=
A_n(r,\psi+c_{10}+c''_5) \subset \wt{U_0}$, by Proposition
\ref{prop:quasiindep}. As claimed after the statement of Proposition
\ref{prop:estimatemeasure}, the series $\sum_{k\in\NN}\wt m_F(A_k)$
diverges. Hence by the above Borel-Cantelli argument, $\limsup_{k}
A_k$ has positive measure.  Since $A_n(r,\psi+c_{10}+c''_5)\subset
A_n(r,\psi +c''_5)$ and by the first claim of Proposition
\ref{prop:geometrictranslation}, the set $\wt E(\psi)$ has positive
$\wt m_F$-measure. Since it is invariant under the geodesic flow and under
$\Ga$, and by ergodicity of the Gibbs measure $m_F$, it has full
measure.

Conversely, assume that the integral $\int_0^{+\infty}
e^{\psi(t)(\delta_0- \delta)} \; dt$ converges, which is still true if
a constant is substracted from $\psi$. Then $\psi(t)\geq c'_5$ for $t$
big enough. Let $A_n= A_n(r,\psi-c'_5) \subset \wt{U_0}$. Again by the
assertion following the statement of Proposition
\ref{prop:estimatemeasure}, the series $\sum_{k\in\NN}\wt m_F(A_k)$
converges. By the standard Borel-Cantelli Lemma, $\limsup_{k}
A_k(r,\psi-c'_5)$ has zero $\wt m_F$-measure. By the second claim of
Proposition \ref{prop:geometrictranslation}, the set $\wt E(\psi)\cap
\wt{U_0}$ has zero $\wt m_F$-measure. Up to taking $R_0$ big enough,
this implies that $\wt E(\psi)$ has zero $\wt m_F$-measure. 
\cqfd

\medskip \rem Let us comment on the range of the numerical constant
$\delta-\delta_0$, crucial for the dichotomy in Theorem
\ref{theo:maintechnicalspiralGibbs}, as the potential $F$ varies. We
only consider the case when $C_0$ is the translation axis of a loxodromic
element of $\Ga$, so that by Remark (2) following the statement of
Theorem \ref{theo:maintechnicalspiralGibbs}, we have, with
$\overline{C_0}$ the image of $C_0$ in $M=\Ga\bs\wt \Ga$,
$$
\delta-\delta_0=P(F) - P(F_{|T^1\overline{C_0}})\;.
$$
\bprop
\begin{enumerate}
\item[(1)] The map $F\mapsto P(F) - P(F_{|T^1\overline{C_0}})$ is
  $1$-Lipschitz for the uniform norm on bounded potentials.
\item[(2)] The set of real numbers $P(F) - P(F_{|T^1\overline{C_0}})$, as $\wt F$
  varies in the set of $\Ga$-invariant bounded Hölder functions on
  $T^1\wt M$, is equal to $]0, +\infty[\,$.
\end{enumerate}
\eprop

\dem For the first observation, the $1$-Lipschitz dependence of
$P(F_{|T^1\overline{C_0}})$ on $F$ is immediate by Equation
\eqref{eq:defPzero}, and so it suffices to show the same for $P(F)$.
This is a direct consequence of our definition of the pressure in
Equation \eqref{eq:defpressure}.  More precisely, given $F_1, F_2$ two
bounded $\Ga$-invariant Hölder-continuous functions on $T^1\wt M$, for every
$\epsilon > 0$, we can choose $m_1, m_2 \in \M$ satisfying
$$
h(m_1) + \int  F_1 \,dm_1 \geq P(F_1) - \epsilon \;\;\hbox{ and }\;\;
h(m_2) + \int  F_2 \,dm_2  \geq P(F_2) - \epsilon\;.
$$
Using the definition of pressure again,  we have  that
$$
P(F_1) \geq h(m_2) + \int  F_1 \,dm_2 \;\; \hbox{ and }\;\;
P(F_2) \geq h(m_1) + \int  F_2 \,dm_1\;.
$$
Comparing these four inequalities  gives that  
$$
\int (F_1 - F_2) \,dm_1 \geq P(F_1) - P(F_2) - \epsilon
\;\; \hbox{ and }\;\;
\int (F_2 - F_1) \,dm_2 \geq P(F_2) - P(F_1) - \epsilon\;,
$$
from which we deduce $|P(F_1)-P(F_2)| \leq \|F_1 - F_2\|_\infty +
\epsilon$.  Letting $\epsilon\ra 0$, this proves that $F\mapsto P(F)$
is $1$-Lipschitz.

\medskip For the second observation, first note that $P(F) -
P(F_{|T^1\overline{C_0}})=\delta-\delta_0$ is positive by Lemma
\ref{lem:gap}.  It now suffices to find two potentials $F, F'$ for
which $P(F) - P(F_{|T^1\overline{C_0}})$ can be arbitrarily large and
$P(F') - P({F'}_{|T^1\overline{C_0}})$ can be arbitrarily close
to $0$.

Given any $L >0$ and a second distinct closed geodesic
$\overline{C_1}$ (which exists since $\Ga$ is nonelementary), we can
choose a bounded potential $F$ on $T^1M$ which is constant with values
$L$ and $0$ on $T^1\overline{C_1}$ and $T^1\overline{C_0}$,
respectively.  If $m_{C_1}$ denotes a probability measure supported on
$T^1\overline{C_1}$ and invariant under the geodesic flow, then by the
definition of the pressure, we have that $P(F_{|T^1\overline{C_0}})=0$
and $P(F) \geq h_{m_{C_1}}(g^{1}) + \int F \,dm_{C_1}= L$, as
required.

Finally, given any $\eta > 0$, we want to construct a bounded
potential $F'$ on $T^1M$ satisfying $P(F') -
P({F'}_{|T^1\overline{C_0}}) < \eta$.  For every $\epsilon >0$, we let
$A_0 = \{v \in T^1M \hbox{ : } d(v, T^1\overline{C_0}) < \epsilon\}$.
We choose $\epsilon \in\;]0,\frac{1}{e}[$  small enough, so that $-
(1-\epsilon)\log(1-\epsilon) - 2\,\epsilon\log\epsilon + 4\,\epsilon\log 2
<\eta$.
We choose $K > h_{\rm top}(g^1)/\epsilon$, and we define a bounded
potential $F'$ on $T^1M$ by $F'(v) = - K \min\{d(v, T^1\overline{C_0})
,1\} \leq 0$.  Given any $m \in \mathcal M$, we can consider two
cases:  Either (a) $m(T^1M - A_0) > \epsilon$ or (b) $m(A_0) \geq
1-\epsilon$.  In case (a),  we have that
$$
h(m) + \int F'\,dm \leq h_{\rm top}(g^1) + \max_{v\in T^1M - A_0} F'(v) 
\leq h_{\rm top}(g^1) - K \epsilon <  0\,.
$$ 
In case (b), we can choose a measurable partition $\alpha =
\{A_n\}_{n\in\NN}$ of $T^1M$, such that:

\smallskip
$\bullet$~ $\alpha$ is {\it generating}, that is, the Borel
$\sigma$-algebra is the smallest $\sigma$-algebra containing
$g^{t_1}A_{i_1} \cap \cdots \cap g^{t_k}A_{i_k}$, for all $k,i_1,
\cdots, i_k \in \NN$ and $t_1, \cdots, t_k \in \RR$;

\smallskip
$\bullet$~ for $n\geq 1$, we have $m(A_n) \leq \epsilon/2^{n-1}$ (note
that $m\big(\bigcup_{n=1}^{+\infty}A_n\big)=1-m(A_0)\leq \epsilon$).

\smallskip
\noindent
If $M$ were compact, then a sufficient condition for the partition to
be generating would be that each element $A_n$, for $n \geq 1$, has
diameter smaller than the injectivity radius of $M$.  (At the level of
the geodesic flow, this is related to choosing the diameter smaller
than the expansivity constant).  More generally, we can assume that
each $A_n$ is the union of suitably separated components, each of which
has diameter smaller than the injectivity radius of points in that
component.  In particular, with $H_m(\alpha)$ the entropy of the
partition $\alpha$ with respect to $m$, we can then bound
\begin{align*}
h(m) + \int F'\,d\mu \leq h(m)  \leq H_m(\alpha) 
\leq -m(A_0)\log m(A_0) - \sum_{n=1}^{+\infty} m(A_n)\log m(A_n)\\
\leq - (1-\epsilon)\log(1-\epsilon) - \sum_{n=1}^{+\infty}
\frac{\epsilon}{2^{n-1}}\log\frac{\epsilon}{2^{n-1}} 
= - (1-\epsilon)\log(1-\epsilon) - 2\,\epsilon\log\epsilon  +
4\,\epsilon\log 2 < \eta\;. 
\end{align*}
In either case, we have that $h(m) + \int F'\;dm < \eta$ and from the
definition, $P(F') - P({F'}_{|T^1\overline{C_0}})=P(F') < \eta$, as
required.

\bigskip Let us now give the main corollary from Theorem
\ref{theo:maintechnicalspiralGibbs}, our logarithm law for Gibbs
measures.

Define the penetration map $\wt \ppp:T^1\wt M\times \RR\ra[0,+\infty]$
of the geodesic lines inside $\Ga \N_{\epsilon_0}C_0$ by $\wt
\ppp(v,t)=0$ if $\pi(\phi_tv)\notin \Ga\N_{\epsilon_0}C_0$, and
otherwise $\wt \ppp(v,t)$ is the maximal length of an interval $I$ in
$\RR$ containing $t$ such that there exists $\ga\in\Ga$ with
$\pi(\phi_sv)\in \ga\N_{\epsilon_0}C_0$ for every $s\in I$. The next
result implies Corollary \ref{coro:mainintro} using Remark (2)
following Theorem \ref{theo:maintechnicalspiralGibbs}.

\bcoro 
Under the assumptions of Theorem \ref{theo:maintechnicalspiralGibbs},
for $\wt m_F$-almost every $v\in T^1\wt M$, we have
$$ 
\limsup_{t\ra+\infty}\frac{\wt \ppp(v,t)}{\ln t} =
\frac{1}{\delta-\delta_0}\;.
$$
\ecoro 

\dem The proof is a standard deduction from Theorem
\ref{theo:maintechnicalspiralGibbs} using the Lipschitz functions
$\psi_n:t\mapsto \kappa \ln(1+t)$ for $\kappa=
\frac{1}{\delta-\delta_0}\pm \frac{1}{n}$, see for instance the proof
of \cite[Theo.~5.6]{HerPau10}. \cqfd

\medskip We end this section by giving a corollary of Theorem
\ref{theo:maintechnicalspiralGibbs} in the special case when $\wt M$
has constant sectional curvature, in a form which is suitable for the
arithmetic applications in the next section. We will use the upper
halfspace model of the real hyperbolic $n$-space $\HH^n_\RR$, whose
boundary at infinity is $\partial_\infty\HH^n_\RR= \RR^{n-1} \cup
\{\infty\}$, and we endow $\RR^{n-1}$ with the usual Euclidean norm
$\|\cdot\|$ and its associated distance. We denote by $x_0$ the point
$(0,1)\in\RR^{n-1} \times\; ]0,+\infty[$. If $\alpha$ is a fixed point
of a hyperbolic element $\ga$ of a given discrete group of isometries
of $\HH^n_\RR$, we denote by $\alpha^\sigma$ its other fixed point,
which does not depend on $\ga$.

\bcoro \label{coro:constcurv} Let $\Ga$ be a nonelementary discrete
group of isometries of $\HH^n_\RR$. Let $\wt F:T^1\HH^n_\RR\ra \RR$ be
a $\Ga$-invariant Hölder-continuous map, with
$\delta=\delta_{\Ga,\,F}$ and $m_F$ finite.  Let $\ga_0$ be a
hyperbolic element of $\Ga$, let $\Ga_0$ be the stabiliser in $\Ga$ of
its translation axis, let $F_0:\Ga_0\bs\wt M\ra \RR$ be the map
induced by $\wt F$, and let $\delta_0=\delta_{\Ga_0,\,F_0}$. Let
$\R_{\ga_0}$ be the set of fixed points in $\RR^{n-1} \cup \{\infty\}$
of the conjugates in $\Ga$ of $\ga_0$. Let $\phi:\;]0,1]\; \ra\;]0,1]$
be a measurable map, such that there exist $c'_2, c'_3 \in\;]0,1[$
such that for every $s,t\in\;]0,c'_2]$, if $s\geq c'_2 \,t$, then
$\phi(s)\geq c'_3\,\phi(t)$.  If $\int_0^1 \phi^{\delta-\delta_0}(s)
/s\,ds$ diverges (respectively converges), then $\mu_{x_0}^F$-almost
every (respectively no) point in $\RR^{n-1}$ belongs to infinitely
many Euclidean balls of centre $\alpha$ and radius $\|\alpha
-\alpha^\sigma \| \phi(\|\alpha -\alpha^\sigma\|)$, as $\alpha$ ranges
over $\R_{\ga_0}$.  \ecoro

\dem Recall that the hyperbolic distance between the horizontal
horosphere at Euclidean height $1$ in $\HH^n_\RR$ and a disjoint
geodesic line with endpoints $x$ and $y$ is $-\ln\frac{\|x-y\|}{2}$,
by a standard hyperbolic distance computation. By the triangle
inequality and the discreteness of $\Ga$, for every compact subset $K$
of $\RR^{n-1}$, there exists $c>0$ such that for every
$\alpha\in\R_{\ga_0}\cap K$ except finitely many of them, we have
$\|\alpha-\alpha^\sigma\|\leq 1$ and, with $C_\alpha$ the geodesic
line with endpoints $\alpha,\alpha^\sigma$,
\begin{equation}\label{eq:heightgeo}
\Big|d(x_0,C_\alpha)-\big|\ln\frac{\|\alpha-\alpha^\sigma\|}{2}
\big|\Big|\leq c\;.
\end{equation}
Let $\psi:t\mapsto -\ln \phi(e^{-t})$ which is a map from
$[0,+\infty[$ to $[0,+\infty[$ satisfying the assumption of the
beginning of Section \ref{sec:spiral} (with $c_2=-\ln c'_2>0$ and
$c_3=-\ln c'_3>0$).

As in \cite[Lem.~5.2]{HerPau10} (and since the Hamenstädt distance on
$\partial_\infty\HH^n_\RR- \{\infty\}= \RR^{n-1}$ is a multiple of the
Euclidean distance), there exists a constant $c'\geq 1$ such that for
every $v\in T^1\wt M$ such that $v_+\in K-(\R_{\ga_0}\cap K)$, we have

$\bullet$~ if $v$ if $(\epsilon_0,\psi)$-Liouville for $(\Ga,\Ga_0)$,
then $v_+$ belongs to infinitely many balls of centre $\alpha$ and
radius $c'\, e^{-d(x_0,\,C_\alpha) -\psi(d(x_0,\,C_\alpha))}$, as
$\alpha$ ranges over $\R_{\ga_0}$.

$\bullet$~ if $v_+$ belongs to infinitely many balls of centre
$\alpha$ and radius $\frac{1}{c'}\, e^{-d(x_0,\,C_\alpha) -\psi(d(x_0,
  \,C_\alpha))}$, as $\alpha$ ranges over $\R_{\ga_0}$, then $v$ if
$(\epsilon_0,\psi)$-Liouville for $(\Ga,\Ga_0)$.

By Equation \eqref{eq:heightgeo}, there exists $c''\geq 1$ such that,
for every $\alpha\in\R_{\ga_0}\cap K$,
$$
\frac{1}{c''}\,\|\alpha -\alpha^\sigma\|\,\phi(\|\alpha -\alpha^\sigma\|)\leq 
e^{-d(x_0,\,C_\alpha) -\psi(d(x_0,\,C_\alpha))}\leq
c''\,\|\alpha -\alpha^\sigma\|\,\phi(\|\alpha -\alpha^\sigma\|)\;.
$$
Since $\int_0^{+\infty}e^{\psi(t)(\delta_0-\delta)}\,dt=\int_0^1
\phi^{\delta-\delta_0}(s)/s\,ds$, the result
follows from Theorem \ref{theo:maintechnicalspiralGibbs}, whose
hypotheses on sum asymptotics are satisfied by the first remark
following its statement (since the curvature of $\wt M$ is constant).
\cqfd

\medskip
\rem As in \cite{ParPau11MZ}, replacing $\HH^n_\RR$ by the Siegel
domain model of the complex hyperbolic space $\HH^n_\CC$, replacing
$\RR^{n-1}$ endowed with the Euclidean distance $\|x-y\|$ by the Heisenberg
group endowed with the Cygan distance $d_{\rm Cyg}(x,y)$, 
the same result holds.

%
%

\section{Arithmetic applications}
\label{sec:arithapplications} 

Let $K$ be either the field $\QQ$ or an imaginary quadratic extension
of $\QQ$, and correspondingly, let $\wh K$ be either $\RR$ or
$\CC$. Let $\OOO_K$ be the ring of integers of $K$. By {\it quadratic
  irrational}, we mean an element in $\widehat K$ which is quadratic
irrational over $K$. For every quadratic irrational $\alpha\in\widehat
K$, let $\alpha^\sigma$ be its Galois conjugate over $K$.

The group $\PSL_2(\wh K)$ acts on $\PP^1(K)=\wh K\cup\{\infty\}$
by homographies, and its subgroup $\PSL_2(\OOO_K)$ preserves the
set $K$ and the set of quadratic irrationals. Though it acts
transitively on the former set, it does not act transitively on the
latter one. Note that, for every quadratic irrational $\alpha$ and
every $\ga\in \PSL_2(\OOO_K)$, we have $(\ga\cdot\alpha)^\sigma=
\ga\cdot (\alpha^\sigma)$. 

Let us fix a finite index subgroup $\Ga$ of $\PSL_2(\OOO_K)$, for
instance a congruence subgroup. We are interested in the approximation
of elements of $\wh K$ by elements in the orbit under $\Ga$ of a fixed
quadratic irrational and of its Galois conjugate.

For every quadratic irrational $\alpha\in\widehat K$, let
$\E_{\alpha,\,\Ga}$ be the (countable, dense in $\widehat K$) set
$\Ga\cdot \{\alpha,\alpha^\sigma\}$, endowed with its Fr\'echet
filter, and let
$$
h(\alpha)=\frac{2}{|\alpha-\alpha^\sigma|}\;.
$$ 
We refer to \cite[\S 6.1]{ParPau11MZ} and \cite[\S 4.1]{ParPau12JMD}
for motivations on this complexity $h(\alpha)$ of a quadratic
irrational $\alpha$, as well as for other algebraic expressions and
comparisons to other algebraic heights. For instance, if $K=\QQ$,
$\Ga= \PSL_2(\ZZ)$ and $\alpha$ is the Golden Ratio
$\frac{1+\sqrt{5}}{2}$, then $\E_{\alpha,\,\Ga}$ is the set of real
numbers whose continued fraction expansion ends with an infinite
string of $1$'s.

Recall that a map $f:[0,+\infty[\;\ra\;]0,+\infty[$ is {\it slowly
  varying} if it is measurable and if there exist constants $B>0$ and
$A\geq 1$ such that for every $x,y$ in $\RR_+$, if $|x-y|\leq B$, then
$f(y)\leq A\,f(x)$. Recall that this implies that $f$ is locally
bounded, hence it is locally integrable; also, if $\log f$ is
Lipschitz, then $f$ is slowly varying.

\btheo\label{theo:diophantine} Let $\alpha_0\in \widehat K$ be a fixed
quadratic irrrational and let $\ga_0\in \Ga$ be a primitive element of
$\Ga$ fixing $\alpha_0$ with $|\ga_0'(\alpha_0)|>1$. Let $\mu^F_{x_0}$
be a Patterson measure on $\wh K\cup\{\infty\}$ associated to a
potential $\wt F$ for $\Ga$ such that $\delta=\delta_{\Ga,\, F}$ and
$m_F$ are finite. Let $\delta_0$ be the critical exponent of
$\ga_0^\ZZ$ for $\wt F$. Let $\varphi:[0,+\infty[\;\ra\;]0,+\infty[$ be
a map such that $t\mapsto \varphi(e^t)$ is slowly varying. If the
integral $\int_{1}^{+\infty} \varphi(t)^{\delta-\delta_0}/t\;dt$
diverges (resp.~converges), then for $\mu^F_{x_0}$-almost every
$x\in\widehat K $,
$$
\liminf_{r\in\E_{\alpha_0,\Ga}} \;
\frac{h(r)}{\varphi(h(r))}\;|x-r|=0\;
({\rm resp.} =+\infty)\;.
$$
\etheo

When $F=0$, this result is due to \cite[Theo.~6.4 (4)]{ParPau11MZ}.

\medskip 
\dem 
Let us first give some details on the notation of this
theorem. Recall (see for instance \cite[Lem.~6.2]{ParPau11MZ}), that
the quadratic irrationals in $\widehat K$ are exactly the fixed points
of the loxodromic elements of $\PSL_2(\OOO_K)$, hence of $\Ga$, since
$\Ga$ has finite index in $\PSL_2(\OOO_K)$. Hence an element $\ga_0$
as in the statement exists and is unique, it is the unique primitive
loxodromic element of $\Ga$ with attractive fixed point $\alpha_0$.

Let $\wt M$ be the real hyperbolic plane $\HH^2_\RR$ if $\wh K=\RR$
and the real hyperbolic space $\HH^3_\RR$ if $\wh K=\CC$. We fix a
point $x_0$ in $\wt M$. Note that $\partial_\infty\wt M=\wh
K\cup\{\infty\}$, and $\Ga$ is a discrete group of isometries
(actually an arithmetic lattice) of $\wt M$, so that a $\Ga$-invariant
potential $\wt F$ on $T^1\wt M$ with $\delta=\delta_{\Ga,\, F}$ does
define a Patterson measure $\mu^F_{x_0}$ seen from $x_0$ (unique up to
scalar multiple if $m_F$ is finite) on $\wh K\cup\{\infty\}$, see
Section \ref{sec:reminder}. Let $\Ga_0$ be the stabiliser of
$\{\alpha_0, {\alpha_0}^\sigma\}$ in $\Ga$ (that is of the translation
axis of $\ga_0$), and let $F_0:\Ga_0\bs T^1\wt M\ra\RR$ be the map
induced by $\wt F$. Since $\ga_0^\ZZ$ has finite index in $\Ga_0$, the
critical exponent $\delta_0$ is equal to $\delta_{\Ga_0,\,F_0}$. Note
that $\E_{\alpha_0,\Ga}$ is exactly the set of fixed points of the
conjugates of $\ga_0$ in $\Ga$.

We may assume that $\varphi\leq 1$. Define $\phi:s\mapsto
\varphi(\frac{2}{s})$, which is a measurable map from
$]0,1]$ to $]0,1[$. The result then follows from Corollary
\ref{coro:constcurv}.  \cqfd

\medskip To conclude, let us give a proof of the last statement of the
Introduction.

\medskip \noindent {\bf Proof of Corollary \ref{coro:approxintro}. }
It is well known that $\Ga_{a,b}$ is a uniform lattice in
$\operatorname{SL}_2(\RR)$ (see for instance \cite[\S 5.2]{Katok92} or
\cite[\S 8.5]{BenPau03}): it is a Fuschian group derived from the
quaternion algebra $\big(\frac{a,\,b}{\QQ}\big)$ over $\QQ$, which is
a division algebra by the nonexistence of nonzero integer solutions to
$x^2-a\,y^2-b\,z^2=0$, hence to $x^2-a\,y^2-b\,z^2+ab \,t^2=0$ by
\cite[Lem.~8.17]{BenPau03}. Let $\Ga= \overline{\Ga_{a,b}}$ be the
image of $\Ga_{a,b}$ in $\PSL_2(\RR)$, which is a cocompact group of
isometries of $\wt M= \HH^2_\RR$, whose action on $\partial_\infty
\HH^2_\RR =\PP_1(\RR)$ is the action by homographies. If $\ga_0=
\begin{pmatrix} x+y\sqrt{a} & z-t\sqrt{a}\\
  b(z+t\sqrt{a}) & x-y\sqrt{a}
\end{pmatrix}$ with $(x,y,z,t)\in\ZZ^4$, then $\operatorname{tr} \ga_0
= 2x$. Hence $|\operatorname{tr}\ga_0|>2$ by the assumptions, that is,
the image of $\ga_0$ in $\Ga$, that we again denote by $\ga_0$, is
hyperbolic. It is well known that its translation length $\ell(\ga_0)$
satisfies (see for instance \cite[page 173]{Beardon83}
$$
\cosh\frac{\ell(\ga_0)}{2}=\frac{|\operatorname{tr}\ga_0\,|}{2}\;.
$$
Let us fix $x_0\in \HH^2_\RR$. By Proposition \ref{prop:hamled}, let
$\wt F:T^1\wt M\ra \RR$ be a $\Ga$-invariant Hölder-continuous map
such that $\mu$ and $\mu^F_{x_0}$ have the same measure class. Since
the conclusion of Corollary \ref{coro:approxintro} depends only on the
measure class of $\mu$, and since $\frac{d(\ga^{-1})_*\mu^F_{x_0}}
{d\,\mu^F_{x_0}}(\ga^+)= \frac{d(\ga^{-1})_*\mu}{d\,\mu}(\ga^+)$ for
every hyperbolic element $\ga\in\Ga$ by
\cite[Théo.~1.c]{Ledrappier95}, as seen in the proof of Proposition
\ref{prop:hamled}, we may assume that $\mu=\mu^F_{x_0}$. Since $\Ga$
is cocompact, both $\delta=\delta_{\Ga,\,F}$ and $m_F$ are finite. Let
$F_0:\Ga_0\bs T^1\wt M\ra \RR$ be the map induced by $\wt F$, and let
$\delta_0= \delta_{\Ga_0,\,F_0}$. By Remark (2) following the
statement of Theorem \ref {theo:maintechnicalspiralGibbs} and by
Equation \eqref{eq:periodcocycle}, we have
$$
\delta_0=\frac{\max\big\{\operatorname{Per}_F(\ga),\,
\operatorname{Per}_F(\ga^{-1})\big\}}{\ell(\ga_0)}=
  \frac1{2\operatorname{arcosh}(\frac{|\operatorname{tr}\ga_0|}{2})}\,
  \max\big\{\frac{d(\ga_0^{-1})_*\mu}{d\,\mu}(\ga_0^+),\,
\frac{d(\ga_0)_*\mu}{d\,\mu}(\ga_0^-)\big\}\;.
$$
Since $\Ga$ is cocompact, by Bowen's period counting theorem (see for
instance \cite[Coro.~1.7]{PauPolSha}), again by Equation
\eqref{eq:periodcocycle}, and by the change of variable
$s=2\cosh\frac{t}{2}$, we have
\begin{align*}
\delta&=\delta_{\Ga,\,F}=\lim_{t\ra+\infty}\; \frac1t
\;\ln\sum_{\ga\in\Ga\;:\;0<\ell(\ga)\leq t}
e^{\operatorname{Per}_F(\ga)} \\ &=\lim_{t\ra+\infty}\; \frac1t
\;\ln\sum_{\ga\in\Ga_{a,\,b},\;\operatorname{tr}(\ga)\neq 0,\pm 2,\;\;
  2\operatorname{arcosh}(\frac{|\operatorname{tr}\ga|}{2})\leq t}
\frac{d(\ga^{-1})_*\mu}{d\,\mu}(\ga^+)\\&=\lim_{s\ra+\infty}\; \frac{1}{2\ln s}
\;\ln\sum_{\ga\in\Ga_{a,\,b},\;2<|\operatorname{tr}(\ga)|\leq s}
\frac{d(\ga^{-1})_*\mu}{d\,\mu}(\ga^+)\;.
\end{align*}
For all $s\geq 0$ and $\epsilon>0$, let $\phi:\;]0,1]\; \ra\;]0,1]$ be
the map $t\mapsto \min\{1,\,\epsilon\, (-\ln t)^{-s}\}$, so that
$\int_0^1 \phi^{\delta-\delta_0}(t) /t\;dt$ diverges if and only if
$s\leq \frac1{\delta-\delta_0}$.

By Corollary \ref{coro:constcurv}, we hence have that if $s\leq
\frac1{\delta-\delta_0}$ (resp.~$s > \frac1{\delta-\delta_0}$),
then, for $\mu$-almost every $x\in\RR$,
$$
\;\liminf_{\alpha\in\Ga_{a,\,b}\cdot \{\ga_0^-,\ga_0^+\}\;:\;
|\alpha-\alpha^\sigma|\ra 0}\;\;
\frac{|x-\alpha|}{|\alpha-\alpha^\sigma|(-\ln|\alpha-\alpha^\sigma|)^{-s}}
\;\leq \;\frac{1}{\epsilon}
\;\;({\rm resp.} \geq \frac{1}{\epsilon})\;.
$$ 
By taking $\epsilon=k$ (resp.~$\epsilon=\frac{1}{k}$) for $k\in\NN$
tending to $+\infty$, this proves the result.  
\cqfd

\bibliography{../biblio}

\bigskip
\bigskip
\noindent {\small 
\begin{tabular}{l} 
D\'epartement de Math\'ematique, Bât.~425, UMR 8628 CNRS\\ 
Université Paris-Sud\\ 
91505 ORSAY Cedex, FRANCE\\ 
{\it e-mail: frederic.paulin@math.u-psud.fr} 
\end{tabular} 
\\ 
 \medskip
\\ 
\begin{tabular}{l}  
Mathematics Institute,
\\  University of Warwick,
\\ Coventry, CV4 7AL, UK
\\{\it e-mail: mpollic@maths.warwick.ac.uk} 
\end{tabular} 
}

\end{document}

%% file: macro.tex
\usepackage{amssymb,epsfig,amsmath,mathrsfs,bbm} 
\usepackage[T1]{fontenc}
\usepackage[dvips]{color}
\usepackage{comment}
\usepackage[utf8]{inputenc}
\usepackage[colorlinks]{hyperref}

\newtheorem{defi}{Definition}[section]
\newtheorem{prop}[defi]{Proposition}
\newtheorem{theo}[defi]{Theorem}
\newtheorem{conj}[defi]{Conjecture}
\newtheorem{lemm}[defi]{Lemma}
\newtheorem{coro}[defi]{Corollary}
\newtheorem{rema}[defi]{Remark}
\newtheorem{exem}[defi]{Example}
\newtheorem{exems}[defi]{Examples}

\newcommand{\bdefi}{\begin{defi}}
\newcommand{\edefi}{\end{defi}}
\newcommand{\bprop}{\begin{prop}}
\newcommand{\eprop}{\end{prop}}
\newcommand{\btheo}{\begin{theo}}
\newcommand{\etheo}{\end{theo}}
\newcommand{\blemm}{\begin{lemm}}
\newcommand{\brema}{\begin{rema}}
\newcommand{\erema}{\end{rema}}
\newcommand{\bexer}{\begin{exem}}
\newcommand{\eexer}{\end{exem}}
\newcommand{\bexems}{\begin{exems}}
\newcommand{\eexems}{\end{exems}}
\newcommand{\bconj}{\begin{conj}}
\newcommand{\econj}{\end{conj}}
\newcommand{\elemm}{\end{lemm}}
\newcommand{\bcoro}{\begin{coro}}
\newcommand{\ecoro}{\end{coro}}
\newcommand{\dem}{\noindent{\bf Proof. }}
\newcommand{\rem}{\noindent{\bf Remark. }}

\renewcommand{\cal}{\mathscr}

\newcommand{\R}{{\cal R}}

\newcommand{\M}{{\cal M}}

\newcommand{\N}{{\cal N}}

\newcommand{\E}{{\cal E}}

\renewcommand{\H}{{\cal H}}
\newcommand{\OOO}{{\cal O}}
\newcommand{\C}{{\cal C}}

\newcommand{\maths}[1]{{\mathbb #1}}  

\newcommand{\RR}{\maths{R}}

\newcommand{\NN}{\maths{N}}
\newcommand{\CC}{\maths{C}}
\newcommand{\QQ}{\maths{Q}}

\newcommand{\HH}{\maths{H}}

\newcommand{\ZZ}{\maths{Z}}
\newcommand{\PP}{\maths{P}}

\newcommand{\ppp}{{\mathfrak p}}

\newcommand{\ra}{\rightarrow}
 
\newcommand{\wt}[1]{{\widetilde{#1}}}
\newcommand{\wh}[1]{{\widehat{#1}}}
\newcommand{\ga}{\gamma}
\newcommand{\Ga}{\Gamma}
\newcommand{\bs}{\backslash}

\renewcommand{\ln}{\log}

\newcommand{\cqfd}{\hfill$\Box$}

\newcommand{\flow}[1]{g^{#1}}
\newcommand{\PSL}{\operatorname{PSL}}

%% file: fig_rightinclusion.pstex_t
\begin{picture}(0,0)%
\includegraphics{fig_rightinclusion.pstex}%
\end{picture}%
\setlength{\unitlength}{3812sp}%
\begingroup\makeatletter\ifx\SetFigFont\undefined%
\gdef\SetFigFont#1#2#3#4#5{%
  \reset@font\fontsize{#1}{#2pt}%
  \fontfamily{#3}\fontseries{#4}\fontshape{#5}%
  \selectfont}%
\fi\endgroup%
\begin{picture}(3780,2904)(191,-2233)
\put(1577,-827){\makebox(0,0)[lb]{\smash{{\SetFigFont{11}{13.2}{\rmdefault}{\mddefault}{\updefault}{\color[rgb]{0,0,0}$g_v(t^-_n)$}%
}}}}
\put(1346,-1376){\makebox(0,0)[lb]{\smash{{\SetFigFont{11}{13.2}{\rmdefault}{\mddefault}{\updefault}{\color[rgb]{0,0,0}$\ga_nx_0$}%
}}}}
\put(2431,-1411){\makebox(0,0)[lb]{\smash{{\SetFigFont{11}{13.2}{\rmdefault}{\mddefault}{\updefault}{\color[rgb]{0,0,0}$D_{\ga_n}$}%
}}}}
\put(313,-151){\makebox(0,0)[lb]{\smash{{\SetFigFont{11}{13.2}{\rmdefault}{\mddefault}{\updefault}{\color[rgb]{0,0,0}$\ga_n\alpha_n x_0$}%
}}}}
\put(3401,-716){\makebox(0,0)[lb]{\smash{{\SetFigFont{11}{13.2}{\rmdefault}{\mddefault}{\updefault}{\color[rgb]{0,0,0}$\pi(v)$}%
}}}}
\put(3481,-1046){\makebox(0,0)[lb]{\smash{{\SetFigFont{11}{13.2}{\rmdefault}{\mddefault}{\updefault}{\color[rgb]{0,0,0}$x_0$}%
}}}}
\put(3456,-1496){\makebox(0,0)[lb]{\smash{{\SetFigFont{11}{13.2}{\rmdefault}{\mddefault}{\updefault}{\color[rgb]{0,0,0}$R_0$}%
}}}}
\put(1551,-2081){\makebox(0,0)[lb]{\smash{{\SetFigFont{11}{13.2}{\rmdefault}{\mddefault}{\updefault}{\color[rgb]{0,0,0}$\ga_n\N_{\epsilon_0} C_0$}%
}}}}
\put(1056,-1081){\makebox(0,0)[lb]{\smash{{\SetFigFont{11}{13.2}{\rmdefault}{\mddefault}{\updefault}{\color[rgb]{0,0,0}$p_{\ga_n}$}%
}}}}
\put(206,-1716){\makebox(0,0)[lb]{\smash{{\SetFigFont{11}{13.2}{\rmdefault}{\mddefault}{\updefault}{\color[rgb]{0,0,0}$\ga_n C_0$}%
}}}}
\put(1095,-907){\makebox(0,0)[lb]{\smash{{\SetFigFont{11}{13.2}{\rmdefault}{\mddefault}{\updefault}{\color[rgb]{0,0,0}$q_n$}%
}}}}
\put(3129,-695){\makebox(0,0)[lb]{\smash{{\SetFigFont{11}{13.2}{\rmdefault}{\mddefault}{\updefault}{\color[rgb]{0,0,0}$v$}%
}}}}
\put(799,329){\makebox(0,0)[lb]{\smash{{\SetFigFont{11}{13.2}{\rmdefault}{\mddefault}{\updefault}{\color[rgb]{0,0,0}$g_v(t_n+\psi(t_n))$}%
}}}}
\put(996, 67){\makebox(0,0)[lb]{\smash{{\SetFigFont{11}{13.2}{\rmdefault}{\mddefault}{\updefault}{\color[rgb]{0,0,0}$g_v(t'_n)$}%
}}}}
\put(1371,-539){\makebox(0,0)[lb]{\smash{{\SetFigFont{11}{13.2}{\rmdefault}{\mddefault}{\updefault}{\color[rgb]{0,0,0}$g_v(t_n)$}%
}}}}
\put(763,-514){\makebox(0,0)[lb]{\smash{{\SetFigFont{11}{13.2}{\rmdefault}{\mddefault}{\updefault}{\color[rgb]{0,0,0}$r$}%
}}}}
\end{picture}%

%% file: fig_leftinclusion.pstex_t
\begin{picture}(0,0)%
\includegraphics{fig_leftinclusion.pstex}%
\end{picture}%
\setlength{\unitlength}{3812sp}%
\begingroup\makeatletter\ifx\SetFigFont\undefined%
\gdef\SetFigFont#1#2#3#4#5{%
  \reset@font\fontsize{#1}{#2pt}%
  \fontfamily{#3}\fontseries{#4}\fontshape{#5}%
  \selectfont}%
\fi\endgroup%
\begin{picture}(3886,2941)(85,-2233)
\put(597,-89){\makebox(0,0)[lb]{\smash{{\SetFigFont{11}{13.2}{\rmdefault}{\mddefault}{\updefault}{\color[rgb]{0,0,0}$q'_n$}%
}}}}
\put(1346,-1376){\makebox(0,0)[lb]{\smash{{\SetFigFont{11}{13.2}{\rmdefault}{\mddefault}{\updefault}{\color[rgb]{0,0,0}$\ga_nx_0$}%
}}}}
\put(2431,-1411){\makebox(0,0)[lb]{\smash{{\SetFigFont{11}{13.2}{\rmdefault}{\mddefault}{\updefault}{\color[rgb]{0,0,0}$D_{\ga_n}$}%
}}}}
\put(100,168){\makebox(0,0)[lb]{\smash{{\SetFigFont{11}{13.2}{\rmdefault}{\mddefault}{\updefault}{\color[rgb]{0,0,0}$\ga_n\alpha_n x_0$}%
}}}}
\put(3481,-1046){\makebox(0,0)[lb]{\smash{{\SetFigFont{11}{13.2}{\rmdefault}{\mddefault}{\updefault}{\color[rgb]{0,0,0}$x_0$}%
}}}}
\put(3456,-1496){\makebox(0,0)[lb]{\smash{{\SetFigFont{11}{13.2}{\rmdefault}{\mddefault}{\updefault}{\color[rgb]{0,0,0}$R_0$}%
}}}}
\put(1551,-2081){\makebox(0,0)[lb]{\smash{{\SetFigFont{11}{13.2}{\rmdefault}{\mddefault}{\updefault}{\color[rgb]{0,0,0}$\ga_n\N_{\epsilon_0} C_0$}%
}}}}
\put(1056,-1081){\makebox(0,0)[lb]{\smash{{\SetFigFont{11}{13.2}{\rmdefault}{\mddefault}{\updefault}{\color[rgb]{0,0,0}$p_{\ga_n}$}%
}}}}
\put(206,-1716){\makebox(0,0)[lb]{\smash{{\SetFigFont{11}{13.2}{\rmdefault}{\mddefault}{\updefault}{\color[rgb]{0,0,0}$\ga_n C_0$}%
}}}}
\put(3129,-695){\makebox(0,0)[lb]{\smash{{\SetFigFont{11}{13.2}{\rmdefault}{\mddefault}{\updefault}{\color[rgb]{0,0,0}$v$}%
}}}}
\put(1577,-827){\makebox(0,0)[lb]{\smash{{\SetFigFont{11}{13.2}{\rmdefault}{\mddefault}{\updefault}{\color[rgb]{0,0,0}$g_v(t^-_n)$}%
}}}}
\put(1179,-58){\makebox(0,0)[lb]{\smash{{\SetFigFont{11}{13.2}{\rmdefault}{\mddefault}{\updefault}{\color[rgb]{0,0,0}$g_v(t^+_n)$}%
}}}}
\put(1150,250){\makebox(0,0)[lb]{\smash{{\SetFigFont{11}{13.2}{\rmdefault}{\mddefault}{\updefault}{\color[rgb]{0,0,0}$g_v(\tau_n)$}%
}}}}
\put(787,585){\makebox(0,0)[lb]{\smash{{\SetFigFont{11}{13.2}{\rmdefault}{\mddefault}{\updefault}{\color[rgb]{0,0,0}$r$}%
}}}}
\end{picture}%

%% file: fig_dynamicalball.pstex_t
\begin{picture}(0,0)%
\includegraphics{fig_dynamicalball.pstex}%
\end{picture}%
\setlength{\unitlength}{3631sp}%
\begingroup\makeatletter\ifx\SetFigFont\undefined%
\gdef\SetFigFont#1#2#3#4#5{%
  \reset@font\fontsize{#1}{#2pt}%
  \fontfamily{#3}\fontseries{#4}\fontshape{#5}%
  \selectfont}%
\fi\endgroup%
\begin{picture}(5791,1351)(1343,-4400)
\put(1651,-3643){\makebox(0,0)[lb]{\smash{{\SetFigFont{11}{13.2}{\rmdefault}{\mddefault}{\updefault}{\color[rgb]{0,0,0}$x_0$}%
}}}}
\put(1958,-3203){\makebox(0,0)[lb]{\smash{{\SetFigFont{11}{13.2}{\rmdefault}{\mddefault}{\updefault}{\color[rgb]{0,0,0}$v$}%
}}}}
\put(1763,-4186){\makebox(0,0)[lb]{\smash{{\SetFigFont{11}{13.2}{\rmdefault}{\mddefault}{\updefault}{\color[rgb]{0,0,0}$R_0$}%
}}}}
\put(2180,-3829){\makebox(0,0)[lb]{\smash{{\SetFigFont{11}{13.2}{\rmdefault}{\mddefault}{\updefault}{\color[rgb]{0,0,0}$v_{\ga\alpha}$}%
}}}}
\put(6744,-4091){\makebox(0,0)[lb]{\smash{{\SetFigFont{11}{13.2}{\rmdefault}{\mddefault}{\updefault}{\color[rgb]{0,0,0}$r$}%
}}}}
\put(6661,-3633){\makebox(0,0)[lb]{\smash{{\SetFigFont{11}{13.2}{\rmdefault}{\mddefault}{\updefault}{\color[rgb]{0,0,0}$\ga\alpha x_0$}%
}}}}
\put(4126,-4336){\makebox(0,0)[lb]{\smash{{\SetFigFont{11}{13.2}{\rmdefault}{\mddefault}{\updefault}{\color[rgb]{0,0,0}$\ga x_0$}%
}}}}
\put(4388,-3892){\makebox(0,0)[lb]{\smash{{\SetFigFont{11}{13.2}{\rmdefault}{\mddefault}{\updefault}{\color[rgb]{0,0,0}$\geq \pi/2$}%
}}}}
\put(4595,-4133){\makebox(0,0)[lb]{\smash{{\SetFigFont{11}{13.2}{\rmdefault}{\mddefault}{\updefault}{\color[rgb]{0,0,0}$p_\ga$}%
}}}}
\end{picture}%

%% file: fig_bishadows.pstex_t
\begin{picture}(0,0)%
\includegraphics{fig_bishadows.pstex}%
\end{picture}%
\setlength{\unitlength}{3812sp}%
\begingroup\makeatletter\ifx\SetFigFont\undefined%
\gdef\SetFigFont#1#2#3#4#5{%
  \reset@font\fontsize{#1}{#2pt}%
  \fontfamily{#3}\fontseries{#4}\fontshape{#5}%
  \selectfont}%
\fi\endgroup%
\begin{picture}(4897,2751)(204,-4516)
\put(2791,-3481){\makebox(0,0)[lb]{\smash{{\SetFigFont{11}{13.2}{\rmdefault}{\mddefault}{\updefault}{\color[rgb]{0,0,0}$\frac{R_0}{3}$}%
}}}}
\put(2971,-3031){\makebox(0,0)[lb]{\smash{{\SetFigFont{11}{13.2}{\rmdefault}{\mddefault}{\updefault}{\color[rgb]{0,0,0}$x_0$}%
}}}}
\put(5086,-3211){\makebox(0,0)[lb]{\smash{{\SetFigFont{11}{13.2}{\rmdefault}{\mddefault}{\updefault}{\color[rgb]{0,0,0}$\approx$}%
}}}}
\put(4861,-3391){\makebox(0,0)[lb]{\smash{{\SetFigFont{11}{13.2}{\rmdefault}{\mddefault}{\updefault}{\color[rgb]{0,0,0}$A_\beta(r)_+$}%
}}}}
\put(811,-3121){\makebox(0,0)[lb]{\smash{{\SetFigFont{11}{13.2}{\rmdefault}{\mddefault}{\updefault}{\color[rgb]{0,0,0}$\approx$}%
}}}}
\put(4239,-3391){\makebox(0,0)[lb]{\smash{{\SetFigFont{11}{13.2}{\rmdefault}{\mddefault}{\updefault}{\color[rgb]{0,0,0}$r$}%
}}}}
\put(586,-3301){\makebox(0,0)[lb]{\smash{{\SetFigFont{11}{13.2}{\rmdefault}{\mddefault}{\updefault}{\color[rgb]{0,0,0}$B_{R_0,\,\beta}$}%
}}}}
\put(219,-2896){\makebox(0,0)[lb]{\smash{{\SetFigFont{11}{13.2}{\rmdefault}{\mddefault}{\updefault}{\color[rgb]{0,0,0}$ \OOO_{\beta x_0}B(x_0,R_0/3)$}%
}}}}
\put(4096,-3031){\makebox(0,0)[lb]{\smash{{\SetFigFont{11}{13.2}{\rmdefault}{\mddefault}{\updefault}{\color[rgb]{0,0,0}$\beta x_0$}%
}}}}
\put(4689,-3031){\makebox(0,0)[lb]{\smash{{\SetFigFont{11}{13.2}{\rmdefault}{\mddefault}{\updefault}{\color[rgb]{0,0,0}$\OOO_{x_0}B(\beta x_0,r)$}%
}}}}
\end{picture}%

%% file: fig_minovoldynball.pstex_t
\begin{picture}(0,0)%
\includegraphics{fig_minovoldynball.pstex}%
\end{picture}%
\setlength{\unitlength}{3729sp}%
\begingroup\makeatletter\ifx\SetFigFontNFSS\undefined%
\gdef\SetFigFontNFSS#1#2#3#4#5{%
  \reset@font\fontsize{#1}{#2pt}%
  \fontfamily{#3}\fontseries{#4}\fontshape{#5}%
  \selectfont}%
\fi\endgroup%
\begin{picture}(6333,1850)(732,-1778)
\put(2018,-1476){\makebox(0,0)[lb]{\smash{{\SetFigFontNFSS{11}{13.2}{\rmdefault}{\mddefault}{\updefault}{\color[rgb]{0,0,0}$z$}%
}}}}
\put(4963,-1606){\makebox(0,0)[lb]{\smash{{\SetFigFontNFSS{11}{13.2}{\rmdefault}{\mddefault}{\updefault}{\color[rgb]{0,0,0}$r$}%
}}}}
\put(903,-1596){\makebox(0,0)[lb]{\smash{{\SetFigFontNFSS{11}{13.2}{\rmdefault}{\mddefault}{\updefault}{\color[rgb]{0,0,0}$R_0$}%
}}}}
\put(5104,-643){\makebox(0,0)[lb]{\smash{{\SetFigFontNFSS{11}{13.2}{\rmdefault}{\mddefault}{\updefault}{\color[rgb]{0,0,0}$z'$}%
}}}}
\put(747,-133){\makebox(0,0)[lb]{\smash{{\SetFigFontNFSS{11}{13.2}{\rmdefault}{\mddefault}{\updefault}{\color[rgb]{0,0,0}$v_-$}%
}}}}
\put(1048,-1091){\makebox(0,0)[lb]{\smash{{\SetFigFontNFSS{11}{13.2}{\rmdefault}{\mddefault}{\updefault}{\color[rgb]{0,0,0}$\frac{R_0}{3}$}%
}}}}
\put(1763,-1056){\makebox(0,0)[lb]{\smash{{\SetFigFontNFSS{11}{13.2}{\rmdefault}{\mddefault}{\updefault}{\color[rgb]{0,0,0}$x_0$}%
}}}}
\put(1526,-741){\makebox(0,0)[lb]{\smash{{\SetFigFontNFSS{11}{13.2}{\rmdefault}{\mddefault}{\updefault}{\color[rgb]{0,0,0}$z''$}%
}}}}
\put(2013,-328){\makebox(0,0)[lb]{\smash{{\SetFigFontNFSS{11}{13.2}{\rmdefault}{\mddefault}{\updefault}{\color[rgb]{0,0,0}$p$}%
}}}}
\put(6958,-166){\makebox(0,0)[lb]{\smash{{\SetFigFontNFSS{11}{13.2}{\rmdefault}{\mddefault}{\updefault}{\color[rgb]{0,0,0}$v_+$}%
}}}}
\put(5338,-1351){\makebox(0,0)[lb]{\smash{{\SetFigFontNFSS{11}{13.2}{\rmdefault}{\mddefault}{\updefault}{\color[rgb]{0,0,0}$\frac{r}{2}$}%
}}}}
\put(5293,-951){\makebox(0,0)[lb]{\smash{{\SetFigFontNFSS{11}{13.2}{\rmdefault}{\mddefault}{\updefault}{\color[rgb]{0,0,0}$\beta x_0$}%
}}}}
\put(1600,-425){\makebox(0,0)[lb]{\smash{{\SetFigFontNFSS{11}{13.2}{\rmdefault}{\mddefault}{\updefault}{\color[rgb]{0,0,0}$v$}%
}}}}
\end{picture}%

%% file: fig_brochette.pstex_t
\begin{picture}(0,0)%
\includegraphics{fig_brochette.pstex}%
\end{picture}%
\setlength{\unitlength}{3812sp}%
\begingroup\makeatletter\ifx\SetFigFont\undefined%
\gdef\SetFigFont#1#2#3#4#5{%
  \reset@font\fontsize{#1}{#2pt}%
  \fontfamily{#3}\fontseries{#4}\fontshape{#5}%
  \selectfont}%
\fi\endgroup%
\begin{picture}(5073,1395)(1371,-1734)
\put(5896,-1051){\makebox(0,0)[lb]{\smash{{\SetFigFont{11}{13.2}{\rmdefault}{\mddefault}{\updefault}{\color[rgb]{0,0,0}$p$}%
}}}}
\put(5041,-1096){\makebox(0,0)[lb]{\smash{{\SetFigFont{11}{13.2}{\rmdefault}{\mddefault}{\updefault}{\color[rgb]{0,0,0}$\leq r$}%
}}}}
\put(6121,-646){\makebox(0,0)[lb]{\smash{{\SetFigFont{11}{13.2}{\rmdefault}{\mddefault}{\updefault}{\color[rgb]{0,0,0}$\leq r$}%
}}}}
\put(1386,-591){\makebox(0,0)[lb]{\smash{{\SetFigFont{11}{13.2}{\rmdefault}{\mddefault}{\updefault}{\color[rgb]{0,0,0}$x_0$}%
}}}}
\put(1901,-1046){\makebox(0,0)[lb]{\smash{{\SetFigFont{11}{13.2}{\rmdefault}{\mddefault}{\updefault}{\color[rgb]{0,0,0}$\pi(v)$}%
}}}}
\put(2051,-531){\makebox(0,0)[lb]{\smash{{\SetFigFont{11}{13.2}{\rmdefault}{\mddefault}{\updefault}{\color[rgb]{0,0,0}$\leq R_0$}%
}}}}
\put(4286,-1231){\makebox(0,0)[lb]{\smash{{\SetFigFont{11}{13.2}{\rmdefault}{\mddefault}{\updefault}{\color[rgb]{0,0,0}$\ga'\alpha' x_0$}%
}}}}
\put(4811,-741){\makebox(0,0)[lb]{\smash{{\SetFigFont{11}{13.2}{\rmdefault}{\mddefault}{\updefault}{\color[rgb]{0,0,0}$p'$}%
}}}}
\put(5446,-486){\makebox(0,0)[lb]{\smash{{\SetFigFont{11}{13.2}{\rmdefault}{\mddefault}{\updefault}{\color[rgb]{0,0,0}$\ga\alpha x_0$}%
}}}}
\put(2601,-781){\makebox(0,0)[lb]{\smash{{\SetFigFont{11}{13.2}{\rmdefault}{\mddefault}{\updefault}{\color[rgb]{0,0,0}$v$}%
}}}}
\end{picture}%

%% file: fig_quasiindependent.pstex_t
\begin{picture}(0,0)%
\includegraphics{fig_quasiindependent.pstex}%
\end{picture}%
\setlength{\unitlength}{3812sp}%
\begingroup\makeatletter\ifx\SetFigFont\undefined%
\gdef\SetFigFont#1#2#3#4#5{%
  \reset@font\fontsize{#1}{#2pt}%
  \fontfamily{#3}\fontseries{#4}\fontshape{#5}%
  \selectfont}%
\fi\endgroup%
\begin{picture}(6320,2053)(263,-2552)
\put(5079,-1641){\makebox(0,0)[lb]{\smash{{\SetFigFont{11}{13.2}{\rmdefault}{\mddefault}{\updefault}{\color[rgb]{0,0,0}$x$}%
}}}}
\put(5042,-1379){\makebox(0,0)[lb]{\smash{{\SetFigFont{11}{13.2}{\rmdefault}{\mddefault}{\updefault}{\color[rgb]{0,0,0}$\ga\alpha x_0$}%
}}}}
\put(4760,-1174){\makebox(0,0)[lb]{\smash{{\SetFigFont{11}{13.2}{\rmdefault}{\mddefault}{\updefault}{\color[rgb]{0,0,0}$r$}%
}}}}
\put(5101,-2027){\makebox(0,0)[lb]{\smash{{\SetFigFont{11}{13.2}{\rmdefault}{\mddefault}{\updefault}{\color[rgb]{0,0,0}$y'$}%
}}}}
\put(5597,-2189){\makebox(0,0)[lb]{\smash{{\SetFigFont{11}{13.2}{\rmdefault}{\mddefault}{\updefault}{\color[rgb]{0,0,0}$r$}%
}}}}
\put(5900,-1630){\makebox(0,0)[lb]{\smash{{\SetFigFont{11}{13.2}{\rmdefault}{\mddefault}{\updefault}{\color[rgb]{0,0,0}$x'$}%
}}}}
\put(5903,-2031){\makebox(0,0)[lb]{\smash{{\SetFigFont{11}{13.2}{\rmdefault}{\mddefault}{\updefault}{\color[rgb]{0,0,0}$\ga'\alpha'x_0$}%
}}}}
\put(624,-993){\makebox(0,0)[lb]{\smash{{\SetFigFont{11}{13.2}{\rmdefault}{\mddefault}{\updefault}{\color[rgb]{0,0,0}$R_0$}%
}}}}
\put(575,-1728){\makebox(0,0)[lb]{\smash{{\SetFigFont{11}{13.2}{\rmdefault}{\mddefault}{\updefault}{\color[rgb]{0,0,0}$\pi(v)$}%
}}}}
\put(3594,-1236){\makebox(0,0)[lb]{\smash{{\SetFigFont{11}{13.2}{\rmdefault}{\mddefault}{\updefault}{\color[rgb]{0,0,0}$y$}%
}}}}
\put(1104,-1843){\makebox(0,0)[lb]{\smash{{\SetFigFont{11}{13.2}{\rmdefault}{\mddefault}{\updefault}{\color[rgb]{0,0,0}$v$}%
}}}}
\put(2964,-1814){\makebox(0,0)[lb]{\smash{{\SetFigFont{11}{13.2}{\rmdefault}{\mddefault}{\updefault}{\color[rgb]{0,0,0}$p$}%
}}}}
\put(3268,-2361){\makebox(0,0)[lb]{\smash{{\SetFigFont{11}{13.2}{\rmdefault}{\mddefault}{\updefault}{\color[rgb]{0,0,0}$q'$}%
}}}}
\put(3553,-2193){\makebox(0,0)[lb]{\smash{{\SetFigFont{11}{13.2}{\rmdefault}{\mddefault}{\updefault}{\color[rgb]{0,0,0}$\ga' x_0$}%
}}}}
\put(3365,-1840){\makebox(0,0)[lb]{\smash{{\SetFigFont{11}{13.2}{\rmdefault}{\mddefault}{\updefault}{\color[rgb]{0,0,0}$p'$}%
}}}}
\put(2049,-670){\makebox(0,0)[lb]{\smash{{\SetFigFont{11}{13.2}{\rmdefault}{\mddefault}{\updefault}{\color[rgb]{0,0,0}$\ga C_0$}%
}}}}
\put(2248,-2488){\makebox(0,0)[lb]{\smash{{\SetFigFont{11}{13.2}{\rmdefault}{\mddefault}{\updefault}{\color[rgb]{0,0,0}$\ga'C_0$}%
}}}}
\put(3124,-1036){\makebox(0,0)[lb]{\smash{{\SetFigFont{11}{13.2}{\rmdefault}{\mddefault}{\updefault}{\color[rgb]{0,0,0}$q$}%
}}}}
\put(2986,-891){\makebox(0,0)[lb]{\smash{{\SetFigFont{11}{13.2}{\rmdefault}{\mddefault}{\updefault}{\color[rgb]{0,0,0}$p_\ga$}%
}}}}
\put(2764,-753){\makebox(0,0)[lb]{\smash{{\SetFigFont{11}{13.2}{\rmdefault}{\mddefault}{\updefault}{\color[rgb]{0,0,0}$\ga x_0$}%
}}}}
\put(948,-1273){\makebox(0,0)[lb]{\smash{{\SetFigFont{11}{13.2}{\rmdefault}{\mddefault}{\updefault}{\color[rgb]{0,0,0}$x_0$}%
}}}}
\end{picture}%